\input ./style/arxiv-ba.cfg
\documentclass[ba,linksfromyear,preprint]{imsart}
\makeatletter
   \@ifpackageloaded{natbib}{}{\usepackage{natbib}}
\makeatother

\pubyear{2015}
\volume{10}
\issue{3}
\firstpage{711}
\lastpage{735}
\doi{10.1214/15-BA944}

\startlocaldefs
\newcommand{\enquote}[1]{``#1''}
\endlocaldefs

\begin{document}

\begin{frontmatter}
\title{Optimal Design in Geostatistics under Preferential Sampling\thanksref{T1}}
\runtitle{Optimal Design in Geostatistics under Preferential Sampling}

\relateddois{T1}{Related articles:
DOI:~\relateddoi[ms=BA944A]{Related item:}{10.1214/15-BA944A},
DOI:~\relateddoi[ms=BA944B]{Related item:}{10.1214/15-BA944B},
DOI:~\relateddoi[ms=BA944C]{Related item:}{10.1214/15-BA944C};
rejoinder at DOI:~\relateddoi[ms=BA944REJ]{Related item:}{10.1214/15-BA944REJ}.}


\begin{aug}
\author[a]{\fnms{Gustavo da Silva} \snm{Ferreira}\corref{}\ead[label=e1]{gustavo.ferreira@ibge.gov.br}}
\and
\author[b]{\fnms{Dani} \snm{Gamerman}\ead[label=e2]{dani@im.ufrj.br}}

\runauthor{G. S. Ferreira and D. Gamerman}

\address[a]{National School of Statistical Sciences, Brazilian Institute of Geography and Statistics, Rio de Janeiro, Brazil, \printead{e1}}
\address[b]{Department of Statistical Methods, Mathematical Institute, Federal University of Rio de Janeiro, Rio de Janeiro, Brazil, \printead{e2}}

\end{aug}

\begin{abstract}
This paper analyses the effect of preferential sampling in
Geostatistics when the choice of new sampling locations is the main
interest of the researcher. A~Bayesian criterion based on maximizing
utility functions is used. Simulated studies are presented and
highlight the strong influence of preferential sampling in the
decisions. The computational complexity is faced by treating the new
local sampling locations as a model parameter and the optimal choice is
then made by analysing its posterior distribution. Finally, an
application is presented using rainfall data collected during spring in
Rio de Janeiro. The results showed that the optimal design is
substantially changed under preferential sampling effects. Furthermore,
it was possible to identify other interesting aspects related to
preferential sampling effects in estimation and prediction in
Geostatistics.
\end{abstract}


\end{frontmatter}


\section{Introduction}

\citet{diggle2010geostatistical} presented a novel methodology to
perform inference in the traditional Geostatistical model under
preferential sampling. They assumed that the sample design could be
described by a log-Gaussian Cox process \citep{moller1998log} and
performed maximum likelihood estimation for the model parameters
through simulation. In addition, they have made simulations to evaluate
the effect of preferential sampling on parameter estimation in
Geostatistics, concluding that this was not negligible. In all
simulations, after obtaining unbiased estimates about the model
parameters, the spatial prediction was made by the {\it plug-in}
method, according to the classical approach to perform inference in
Geostatistics \citep[see][]{cressie1993statistics,diggle1998model}.

Based on this work, many others questions have emerged, and one of them
was related to the influence of preferential sampling in the optimal
design choice. This procedure is widespread in Geostatistics
literature. There are several recent papers dealing with this, such as
\citet{zidek2000designing,fernandez2005effect,zhu2005spatial,diggle2006bayesian,%
gumprecht2009designs,boukouvalas2009,muller2010compound}, among others.
The advances made by \citet{muller1999simulation,muller2004optimal} and
\citet{muller2007simulation} that propose methods based on maximization
of utility functions are especially relevant. Since these procedures
incorporate the Bayesian approach and intensive simulation methods in a
natural way, they became quite appropriate to perform the optimal
design choice in Geostatistics.\vadjust{\goodbreak}

Here the uncertainties about model parameters will be considered
allowing us to evaluate their impact on spatial prediction (or kriging)
under preferential sampling. With this aim, spatial prediction about
the underlying process are performed by analysing the predictive
distribution conditional on the observed data and the design sample.
Then, comparisons between these predictions and those obtained
according to the classical approach are made.

However, the most important contribution of this paper is its analysis
of preferential sampling effects in the process of obtaining the
optimal design in Geostatistical models. Using an approach based on
maximizing utility functions \citep{muller1999simulation} to obtain
optimal design, the influence of preferential sampling is evaluated in
situations where the researcher's goal is to optimize some objective
function, e.g. reduce predictive variances. It will be shown through
simulations that the optimal decision about this choice is
substantially modified under a preferential sampling effect.

This paper is organized as follows: Section 2 presents some background
about spatial processes and a methodology for fully Bayesian inference,
and spatial prediction in Geostatistics under preferential sampling
using MCMC methods is described. Section~3 presents a method to obtain
the optimal design based on maximization of utility functions. Section
4 combines this methodology to obtain the optimal design under
preferential sampling. In Sections 5 and 6, we use the methodology to
obtain the optimal design in some simulated examples, and we analyse a
real dataset where the researcher is interested in monitoring the
occurrence of extreme events. Finally, Section 7 discusses the results
obtained.\vspace*{-1pt}

\section{Spatial Process \& Optimal Design}

This section presents some background on Geostatistics and point
processes, and pre\-sents a procedure for obtaining the optimal design
through utility functions.\vspace*{-1pt}

\subsection{Geostatistics}

Geostatistics deals with stochastic processes defined in a region $D$,
$D \subset \Re^k$, where usually $k=$ 1 or 2. Following the approach in
\citet{diggle1998model}, one can assume that the researcher is
interested in studying the features of the stochastic process
\mbox{$\{S(x):x\in D\}$}. Additionally, assume that $S(x)$ is a Gaussian
stationary and isotropic process, with zero mean, constant variance
$\sigma^2$ and autocorrelation function
$\rho(S(x),S(x+h);\phi)=\break\rho(\parallel h \parallel;\phi)$, $\forall x
\in D$, which may depend on one or more parameters represented by
$\phi$. Several families of autocorrelation function can be found in
the literature \citep[see][]{cressie1993statistics,diggle2007model}.

Assuming that $n$ observations $Y_i = Y( x_i )$, $i= 1,\dots,n$ are
available and
$$Y_i=\mu + S(x_i)+Z_i,$$
$$E[Z_i]=0, Var[Z_i]=\tau^2,\quad \forall i,\vadjust{\goodbreak}$$
one can consider ${\bf Y} = ( Y_1, \dots, Y_n )$ as a noisy
version of the underlying process $S$. As usual, it will be assumed
that ${\bf Z} = (Z_1, \dots, Z_n )$ has a Gaussian distribution
independent of the process $S$.

Under the Bayesian paradigm, the posterior distribution of all model
parameters must be obtained to make inference. Assuming that
$\theta=(\tau^2,\sigma^2,\phi,\mu)$ represents the set of unknown
parameters and that $p({\bf y} \mid \theta)$ is the likelihood
function, we need to specify a prior distribution $p(\theta)$ to obtain
the posterior distribution $p(\theta \mid {\bf y})$. It is usual to
assign Gamma distributions to the parameters $\phi$, $\sigma^{-2}$ and
$\tau^{-2}$ and a Gaussian distribution to~$\mu$. Since the posterior
distributions obtained have no closed form, we can approximate them
using MCMC methods.

Usually there is an additional interest in obtaining predictions of the
process $S$ at locations not observed. Computations required in
inference are usually implemented via discretization of the region $D$
in $M$ subregions or cells. Thus, we can redefine the underlying
Gaussian process over the centroids of these subregions,
$S=\{S_1,\dots,S_M\}$. In addition, we define the partition $S=\{S_{\bf
y},S_N\}$ to distinguish the underlying process associated with the $n$
observed and the $N=M-n$ unobserved locations. The distribution of $S$
is a multivariate Gaussian distribution with dimension $M=n+N$, mean
vector {\bf 0} and autocorrelation matrix given by
$$R=\left(\begin{array}{cc} R_{n} & R_{n,N} \\ R_{N,n} & R_{N} \end{array}\right),$$
whose elements are defined by the autocorrelation function
$\rho(\cdotp;\phi)$. Then, we have the model
%
\begin{equation}
[{\bf y} \mid S_{\bf y},\mu,\tau^{-2}] \sim N({\bf 1}\mu + S_{\bf y},\tau^2 I_{n}), \label{obs}
\end{equation}
\begin{equation}
[S \mid \sigma^{-2},\phi] \sim N({\bf 0}, \sigma^2R) \label{intensity}
\end{equation}
that is completed with the priors $\mu \sim N(0,k)$,
$\tau^{-2} \sim G(a_{\tau},b_{\tau})$ and $\sigma^{-2} \sim
G(a_{\sigma},b_{\sigma})$, where $k, a_{\tau},b_{\tau},a_{\sigma}$ and
$b_{\sigma}$  are known hyperparameters.

Without loss of generality with respect to the objectives of this work,
the exponential autocorrelation function $\rho(\parallel h
\parallel;\phi)=\exp (-\parallel h \parallel / \phi )$ will
be used in the sequel.

The full conditional distributions for the parameters $\mu,\tau^{-2}$,
$\sigma^{-2}$, and $S$ can be updated by Gibbs sampling steps, whereas
the range parameter $\phi$ can be updated by Metropolis steps in the
MCMC algorithm. The full conditional distribution for $S$ is of
particular interest in Geostatistics and has Gaussian distribution with
mean vector and covariance matrix given by
$$\sigma^2 R_{N,n}(\tau^2 I_{n} + \sigma^2R_{n})^{-1}({\bf y}-{\bf 1}\mu)$$
and
$$\sigma^2 R_{N}-\sigma^2 R_{N,n}(\tau^2 I_{n} + \sigma^2R_{n})^{-1}\sigma^2 R_{n,N}.$$

Expressions above are known as {\it kriging predictor} and {\it kriging
variance}, respectively. There are extensions of the basic model of
Geostatistics, most of them developed to deal with non-stationarity
(\citet{higdon1999non,fuentes2001new,fuentes2002spectral} and
\citet{bornn2012}) or non-Gaussianity \citep{diggle1998model}.

\subsection{Spatial Point Process}

The use of point processes for modelling patterns of points in space
has been intensified in the last decades, specially after the
publication of the classic texts of \citet{ripley2005spatial} and
\citet{diggle1983}. The study of point processes has also evolved,
based in recent computational advances. \citet{moller2007modern}
present an excellent review of the methods and models to spatial point
processes and highlight several applications and computational aspects
related to inference.

A spatial point process {\bf X} can be understood as a random finite
subset of event locations belonging to a certain limited region $D
\subset \Re ^k$, where usually $k=1,2$ or~3. A~spatial point process,
defined in a region $D \subset \Re^k$, governed by a non-negative
random function $\Lambda=\{\Lambda(x): x \in \Re^k\}$, is a {\it Cox
process} if the conditional distribution of $[{\bf X} \mid
\Lambda=\lambda]$ is a Poisson process with intensity function
$\lambda(x)$. Additionally, if one can assume that $\Lambda(x)=\exp
\{Z(x)\},$ where $Z=\{Z(x) \in \Re^k\}$ is a stationary and isotropic
Gaussian process, then it is said to be {\it log-Gaussian Cox process}
\citep{moller1998log}.

The likelihood function associated with this process is given by
$$p({\bf x} \mid Z) \propto \exp \left(-\int_{D} \exp \{Z(x)\}dx\right) \prod_{i=1}^{n} \exp \{Z(x_i)\}.$$
Although this function is not analytically tractable, inference on a
log-Gaussian Cox process can be performed in a reasonably simple way
through Monte Carlo simulation methods \citep{moller1998log}. Again, it
is usual to represent the domain $D$ of the point process as a grid and
approximate the Gaussian process $Z(x)$ by a finite-dimensional normal
distribution defined on the grid. \citet{waagepetersen2004convergence}
showed that the discretized posterior for log-Gaussian Cox process
converges to the exact posterior when the sizes of the grid cells tend
to zero. According to \citet{moller1998log}, if the point process is
reasonably aggregated and has moderate intensity, the choice of the
grid does not need to be fine to produce good results in inference. The
inference conducted under the Bayesian paradigm, using MCMC methods,
makes it relatively easy to obtain approximations of $p({\bf x} \mid
Z)$. On the other hand, the computational cost may be high
\citep{moller2007modern}.

\subsection{Geostatistics under Preferential Sampling}

In the Geostatistics literature, it is common to consider the sample
points {\bf x} as fixed or, if coming from a stochastic process,
independent of the process $S(x)$. When the sample design is
stochastic, we must specify the joint distribution of $[Y, S, X]$. We
have a process under {\it preferential sampling} if $[S, X] \neq [S]
[X]$, i.e. the sampling design is dependent of the spatial process.

The class of models proposed by \citet{diggle2010geostatistical} to
accommodate preferential sampling effect assumes that, conditional on
$S$, $X$ is a {\it log-Gaussian Cox process} with intensity
$\lambda(x)=\exp\{ \alpha + \beta S(x) \}$. In addition, conditional on
$[S,X]$, we have that $Y_{i} \sim N[\mu+S(x_i),\tau^2]$, $i=1,\dots,n$.
\citet{diggle2010geostatistical} presents a way to evaluate this
distribution through a fine discretization of the region $D$. The
region $D$ can be discretized into $M$ cells with centroids $x_i, i =
1, \dots, M$, where only one point is expected in each cell.

In a Bayesian approach, we need to obtain the posterior distributions
to make inference\ about the model parameters. For this purpose we
assign Gaussian priors to $\alpha$ and $\beta$.
\citet{pati2011bayesian} proved that the use of improper priors for the
parameter that controls the preferential sampling effects produces
proper posteriors. Thus, the data provide enough information to perform
inference for this model parameter, even under vague prior information.

Since the posterior distributions of these parameters have no closed
form, we use an approximation making use of MCMC methods in a
discretized version of the model given by
$$p({\bf x} \mid S,\alpha,\beta) \propto \prod_{i=1}^{M} [\exp(\alpha + \beta S(x_i))]^{n_i}
\exp \left(-\sum_{i} \Delta_i \exp(\alpha+\beta S(x_i))\right),$$
where $n_i$ and $\Delta_i$ represent the counts and the volume of the
$i$th subregion, $i =1, \dots, M$.

The full conditional distributions of the parameters $\mu, \phi,
\tau^{-2}$ and $\sigma^{-2}$ are given by the same expressions obtained
in the case of non-preferential sampling. The full conditional
distributions of $\beta$ and $\alpha$ are updated in a similar way.
The full conditional distribution of $S$ is updated in MCMC by Metropolis
steps (see details in Appendix A).

\section{Optimal Design}

In general, finding an optimal design involves procedures for obtaining
the maximum or minimum of an objective function. These objective
functions usually quantify the gains and losses related to each
possible decision. In this case, we need to decide which locations (in
time or space) will be collected in order to better understand certain
characteristics of the phenomenon. When the phenomenon of interest is
studied assuming that it is governed by an underlying stochastic
process, the methodology for obtaining the optimal design is usually
performed via Decision Theory. For more details related to Decision
Theory, see the classical textbook of \citet{degroot2005optimal}.

According to \citet{muller1999simulation}, the procedure for obtaining
the optimal design can be performed by defining a utility function
$u({\bf d},\theta,{\bf y_d})$, where ${\bf d}=(d_1,\dots,d_m)'$
represents the $m$ new sample locations, $d_i \in D$, and ${\bf
y_d}=(y_{d_1},\dots,y_{d_m})'$ is the vector of future observations
arising from it, $i=1,\dots,m$. After a set of observations {\bf y} is
available, the optimal design is the vector ${\bf d}^{*}$ that
maximizes the function
$$U({\bf d})= \int u({\bf d},\theta,{\bf y_d})p_{{\bf d}}({\bf y_d} \mid \theta,
{\bf y})p(\theta \mid {\bf y})d\theta d{\bf y_d}=E_{\theta,{\bf y_d} \mid
{\bf y}}[u({\bf d},\theta,{\bf y_d})].$$

In other approaches, one can include additional information to obtain
the optimal design. In order to achieve the optimal design in time,
\citet{stroud2001optimal} used covariates to reduce the model parameter
uncertainties, thereby obtaining the more appropriate point in time for
the return of a patient undergoing treatment. On the other hand,
\citet{ding2008bayesian} used a hierarchical model to relate the
effects of different treatments in clinical studies to determine the
optimal design associated with a specific treatment.

An interesting strategy to optimize the expected utility is known as
{\it an augmented model} \citep{muller1999simulation}. In this case,
the optimal design point is considered as a parameter and can be
estimated through its posterior distribution. Thus, an artificial
probability model $h({\bf d},\theta,{\bf y_d})$ is defined assuming
that $D$ is bounded and $u({\bf d},\theta, {\bf y_d})$ is non-negative
and limited. The distribution $h({\bf d},\theta,{\bf y_d})$ is given by
$$h({\bf d},\theta,{\bf y_d}) \propto u({\bf d},\theta,
{\bf y_d})p_{{\bf d}}({\bf y_d} \mid \theta)p(\theta).$$

Under this assumptions, the marginal distribution of ${\bf d}$ is
proportional to
$$h({\bf d}) \propto \int u({\bf d},\theta,{\bf y_d})p_{{\bf d}}({\bf y_d} \mid \theta)
p(\theta)d\theta d{\bf y_d} = U({\bf d}).$$
Therefore, finding the mode of $h({\bf d})$ becomes equivalent to
maximizing the expected utility of ${\bf d}$. Since the sampled values
of ${\bf d}$ in the MCMC algorithm will be concentrated near regions of
high expected utility, less time is consumed in the simulation
procedure. In other words, one only needs to find the mode of the
distribution, instead of dealing with an optimization problem.

It is important to remark that this methodology can be naturally
performed in situations where one wants to choose $m>1$ new design
points. The main difficulty in this case is to face the increased
computational cost involved since the evaluation of $u(\bf d)$ is
required at each iteration in MCMC. In addition, a large sample of
${\bf d}=(d_1,\dots,d_m)'$ would be needed to evaluate $h(\bf d)$
properly. Alternatively, it is also possible to obtain the $m$ new
points sequentially. However, in this case, there is no guarantee that
the resulting design will be optimal. The evaluation of the utility
function also allows removing points instead of adding new points. In
this case, the decision would be based on gains and losses in the
utility after eliminating each site or group of sites.

\subsection{Utility Functions}

Commonly in Geostatistics, the researcher focuses more on spatial
prediction than on inference about its parameters. In these cases, the
choice of a utility function that depends on the predicted values and
their variances is a natural choice. In this sense, the prediction
variance associated with the location $x \in D$, $V(S(x) \mid {\bf y})$
gives to the researcher the degree of uncertainty about his/her
predictions after observing the data {\bf y}. Thus, reducing this
variance throughout the region $D$ constitutes a reasonable criterion
for choosing a new sample location $d$, $d \in D$. Based on this
principle, the following utility function is a natural choice for
implementing the approach described in \citet{muller1999simulation} in
the context of Geostatistics:
%
\begin{equation}\label{eq:utilidade1}
u({\bf d},\theta,{\bf y_d})=\int [V(S(x) \mid \theta,{\bf y})-V(S(x)
\mid \theta, {\bf y}, {\bf y_d})]dx,
\end{equation}
which can be interpreted as the gain obtained in reduction of
prediction uncertainty after observing ${\bf y_d}$. Thus we need to
maximize its expectation given by $U({\bf d})$.

Finally, it is also important to note that the predictive variances
used in the utility function $u({\bf d}, \theta,{\bf y_d})$ depend only
on the location of ${\bf y_d}$ in $D$, instead of its value. This
feature of the predictive variance avoids the necessity of evaluating
all possible values of ${\bf y_d}$ while maximizing $U({\bf d})$, thus
reducing the computational cost involved. More details about the
evaluation of the utility function (\ref{eq:utilidade1}) are presented
in Appendix B.

Other utility functions could be chosen. In Section 6, we employed a
utility function that is higher for locations where extreme values or
exceedances are expected. In this case, we have
%
\begin{equation}\label{eq:utilidade2}
u(d,\theta,y_d)=P[|S(x_d)| > x_0 \mid \theta,{\bf y}]
\end{equation}
where $x_d$ represents the location associated with $y_d$ and $x_0$ is
an extreme threshold. In practical situations, this expression could be
improved by considering costs and risks related to the occurrences.

\section{Optimal Design under Preferential Sampling}

The problem of obtaining the optimal design ${\bf d}^*$ for spatial
processes under preferential sampling can be performed based on
optimization of
$$U({\bf d})=E_{\theta,{\bf y_d} \mid {\bf x},{\bf y}}[u({\bf d},\theta,
{\bf y_d})] = \int u({\bf d},\theta,{\bf y_d})p({\bf y_d} \mid \theta,
{\bf x},{\bf y})p(\theta \mid {\bf x},{\bf y})d\theta d{\bf y_d}$$
where $p(\theta \mid {\bf x},{\bf y})$ is obtained in Section 2. Given
the posterior distribution of $\theta$ and the utility function, we can
achieve the optimal design by seeking the mode of the posterior
pseudo-distribution of ${\bf d}$. The greatest difficulty in this step
is evaluate the effects of preferential sampling in the utility
function $u({\bf d},\theta,{\bf y_d})$. In many cases, like the
situations presented in this paper, the evaluation of this function
becomes challenging.

As will be noted in the next section, preferential sampling directly
impacts on the estimation of the mean $\mu$ of the underlying Gaussian
process. If the utility function $u({\bf d},\theta,{\bf y_d})$ depends
on this parameter, the choice of optimal design will be greatly
affected. On the other hand, it would be expected that preferential
sampling will also affects utility functions defined to quantify
reductions of uncertainty related to each choice. This actually happens
since the spatial configuration of the sample points also provides
information about the underlying process.

Using as an example the utility function defined in Section 3, one
would need to include the information provided by the observed point
process {\bf x}, i.e.
$$u({\bf d},\theta,{\bf y_d})=\int [V(S(x) \mid \theta,{\bf y},{\bf x})-V(S(x)
\mid \theta, {\bf y},{\bf x},{\bf y_d})]dx,$$
and one would still need to know the variance of the distribution of
$[S \mid \theta,{\bf y},{\bf x}]$. Samples of this distribution are
easily obtained during the implementation of MCMC, as described in
Section 2, but we cannot directly obtain estimates of this variance at
each iteration of the algorithm. To deal with this difficulty, one must
resort to an approximation. A sampling-based approximation would
require an additional MCMC sub-chain at each iteration, thus increasing
substantially the already high computational cost. An analytic, cheaper
alternative is to use a Gaussian approximation of this distribution in
order to evaluate the variances required in the utility function
(\ref{eq:utilidade1}). Note that this approximation is only used to
evaluate $u({\bf d},\theta,{\bf y_d})$. More detail about this
approximation are presented in Appendix C. On the other hand, if the
utility function (\ref{eq:utilidade2}) is used, we can directly
evaluate $U(d)$ from $[S \mid {\bf y,x}]$, since a sample of this
distribution will be available after performing the MCMC.

\section{Simulation Study}

In this section, we simulate datasets according the model presented in
Section 2.3 to illustrate the effects of preferential sampling on
inference and optimal design choice. In all cases, we chose a set of
parameters which do not produce very large samples, in order to enhance
the effects of preferential sampling. We also consider situations where
one needs to add $m=1$ and $m=2$ new locations to the sample design. In
order to produce the observations $y_1,\dots,y_n$, we first generate a
surface $S(x)$, according the Geostatistical model presented in Section
2.1, over the discretized region $D$. Conditioning on $S$, we then
generate a point process from the log-Cox Gaussian model presented in
Section~2.2 in order to obtain $n$ observations. We considered five
simulated studies:
{\setlength\leftmargini{4.3em}
\begin{itemize}
\item[Case I.] Simulation considering a one-dimensional region
    $D=[0, 100]$, partitioned into $M=100$ sub-regions, where
    $(\alpha; \beta; \mu; \sigma^2; \phi; \tau^2; n)=(-3; 2; 12; 2;
    20; 0.1; 18)$.

\item[Case II.]  Simulation considering a one-dimensional region
    $D=[0, 100]$, partitioned into $M\,{=}\,200$ sub-regions, where
    $(\alpha; \beta; \mu; \sigma^2; \phi; \tau^2; n)\,{=}\,(-3.5; 3; 12;
    1; 20; 0.01; 123)$.

\item[Case III.]  Simulation considering a one-dimensional region
    $D=[0, 200]$, partitioned\break  into $M=200$ sub-regions, where
    $(\alpha; \beta; \mu; \sigma^2; \phi; \tau^2; n)=(-1.5; 0.5;
    12; 1; 20;\break 0.01;  56)$.

\item[Case IV.]  Simulation considering a two-dimensional region
    $D=[0, 100]$, partitioned into $M=225$ sub-regions, where
    $(\alpha; \beta; \mu; \sigma^2; \phi; \tau^2; n)=(-8; 2; 12; 2;
    20; 0.1; 12)$.

\item[Case V.]  Exactly as Case IV but with $M=400$ sub-regions.
\end{itemize}}
{\it Case V} was considered in order to evaluate the discretization
effect. Figure 1 shows the Gaussian processes $S$ simulated for {\it
Cases I--IV}. Since $\beta> 0$ in all simulations, the observations are
concentrated near the sites where the process $S$ showed higher values.

\begin{figure}[!htb]
\centering
\includegraphics{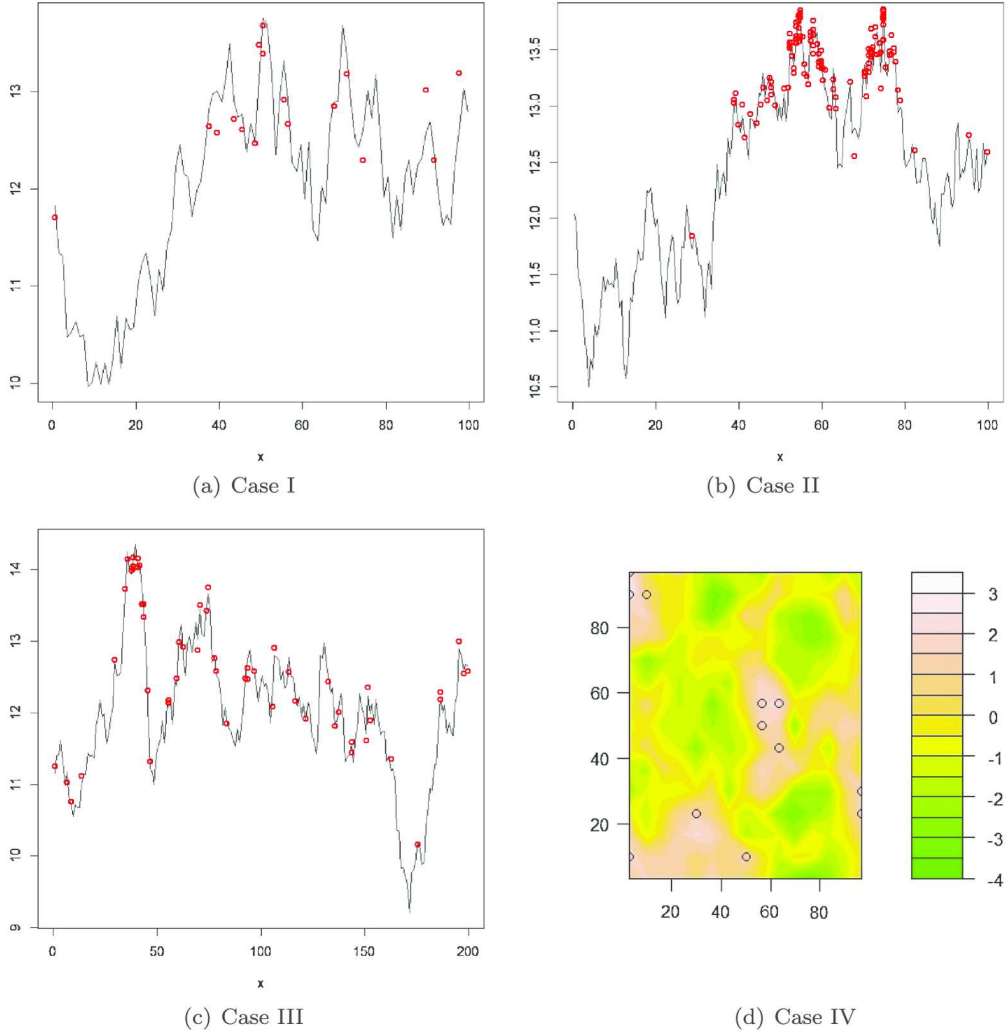}
\caption{Realization of the one-dimensional (Cases I, II and III) and the two-dimensional
(Case IV) simulated processes with respective observed values
(circles). Except for  (c), the observations are concentrated near the
sites where the process S shows higher values.}
\end{figure}

\subsection{Inference and Prediction}

Posterior inference about model parameters was performed. Inference was
also performed  without considering preferential sampling, i.e. using
the Geostatistical model of Section 2.1 to allow comparison. Prior
distributions $\mu \sim \alpha \sim \beta \sim N(0;10^3)$, $\tau^{-2}
\sim \sigma^{-2} \sim G(2; 0.5)$ and $\phi \sim G(2; 0.05)$ were used
in all cases. Furthermore, in {\it Cases I--III}, 500,000 iterations
were generated in the MCMC algorithm and only the last 100,000 were
used to compose the posterior distribution samples of the model
parameters. For the other cases, 400,000 iterations were generated and
only the last 50,000 were considered. The convergence of the chains was
assessed by visual inspection of several chains generated from
different initial values.

Considering the effect of preferential sampling, the posterior
distributions of each model parameters are concentrated around the true
values of the model's parameters, except for the $\sigma^2$ in {\it
Case I}. However, this parameter was also underestimated considering a
non-preferential model. Figure 2 shows that the variogram estimated
assuming the preferential sampling effects is slightly closer to the
true variogram for the {\it Cases~I} and $V$. Figure 2 also shows that
the variograms estimated by the preferential model are closer than
those estimated by the non-preferential model in {\it Cases II} and
$IV$ (Figure~2(a),(b),(e), and (f)). Finally, in {\it Case III}, since
$\beta$ is small, the estimated variograms are very similar (also in
Figure 2).

In particular, in {\it Case IV}, the posterior distributions are not
concentrated around the true values for all parameters of the model.
This difficulty is partly justified by the small sample size. However,
the results obtained by the model with preferential sampling are more
satisfactory. Besides a better estimated variogram, this also can be
observed through the posterior distributions of $\mu$, shown in Figure
3.

Figure 4 shows the predicted values of $S$, represented by the median
of a posterior~$S$, and the respective 95\% credibility intervals for
each model in {\it Cases I} and $II$.

Analysing Figure 4(a)--(b), we can see that only the credible intervals
which consider the effect of preferential sampling encompass the most
extreme values of the simulated process $S$. Additionally, it can be
seen that the point estimates of $S$ in the regions where the process
was hardly observed are better with the model with preferential
sampling. The differences are even more pronounced when we analyse the
results for {\it Case II}, where the intervals are narrower for the
preferential model (Figure 4(c)--(d)). In addition, the
non-preferential model underestimates the process $S$.

Since the estimated variograms obtained by both models were similar in
{\it Cases I} and~$V$, it seems reasonable to conclude that the
differences observed in prediction for these cases were caused mostly
due to the differences between the predictive distributions  $[S\mid
{\bf y}]$ and $[S\mid {\bf y,x}]$.

Finally, Figure 5 shows the predicted surfaces of $S+\mu$, represented
by the posterior means, obtained for each model in {\it Case IV}. It
can be observed that only the preferential model can identify regions
where the underlying process presents low values. This feature avoid
predictions concentrated around the mean, like those obtained by the
traditional kriging methods. In addition, we can conclude that the
small sample size intensifies the effects of preferential sampling in
prediction.
\begin{figure}
\centering
\includegraphics[scale=0.95]{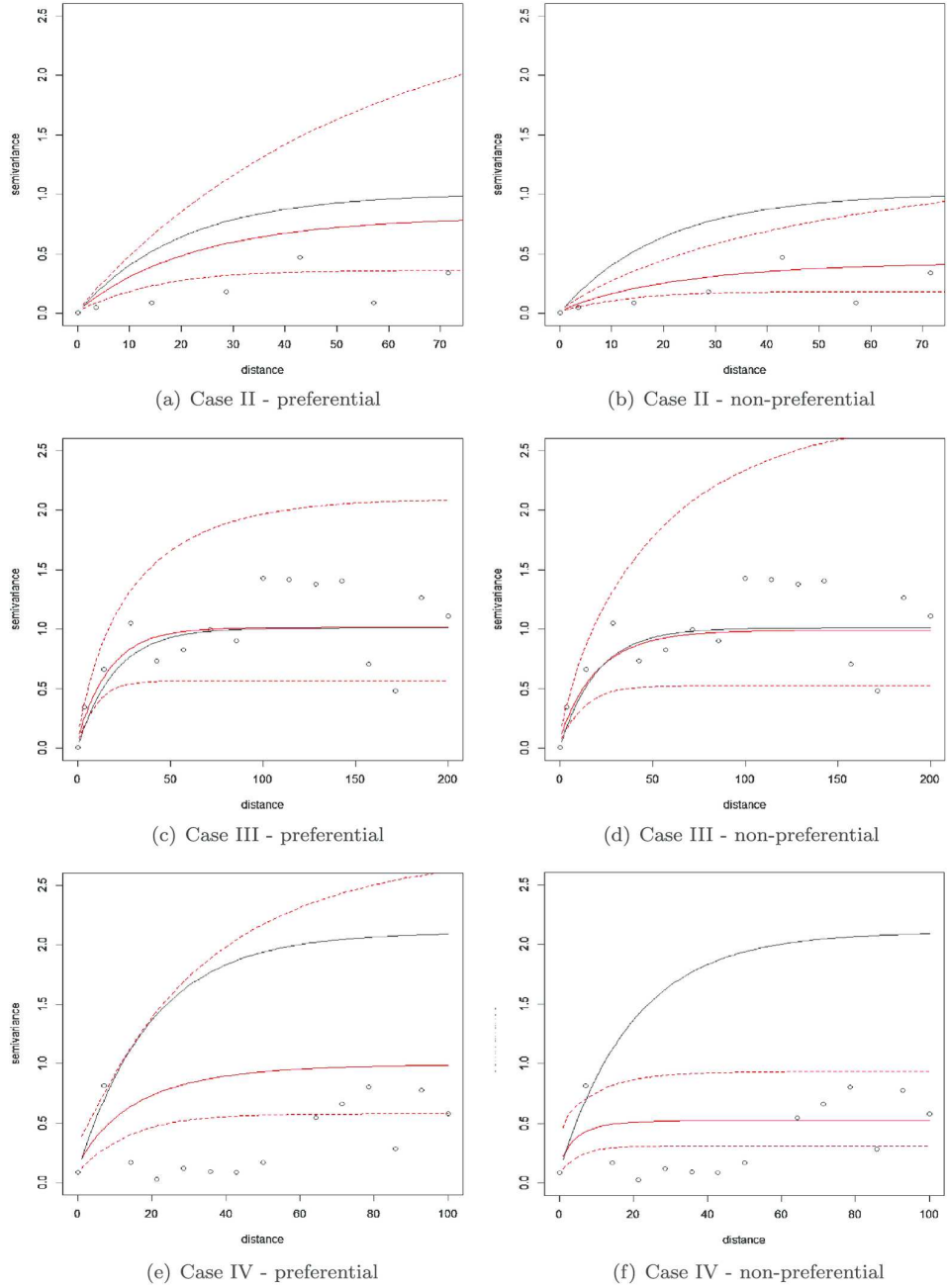}
\caption{Posterior medians (red line) and respective 95\% credibility intervals
(dashed lines) of the variogram obtained by the preferential model
(left) and by the non-preferential model (right) for Cases II, III and
IV. The circles represent the empirical variogram and the black line
represents the true variogram.}
\end{figure}

\begin{figure}[!htb]
\centering
\includegraphics{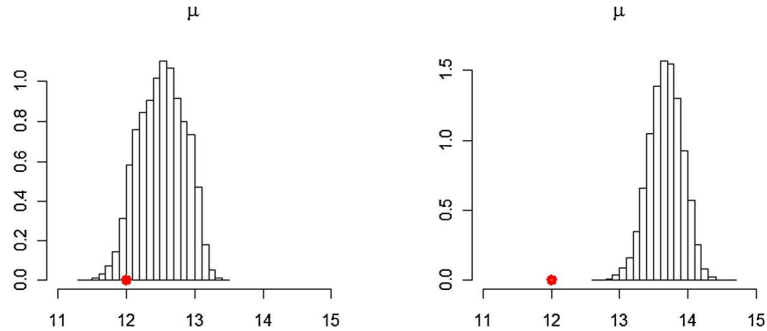}
\caption{Posterior distributions of $\mu$ under preferential sampling
(left) and without considering this effect (right) for {\it Case IV}
(red points represent the true values of the parameters).}
\end{figure}

\begin{figure}[!htb]
\centering
\includegraphics{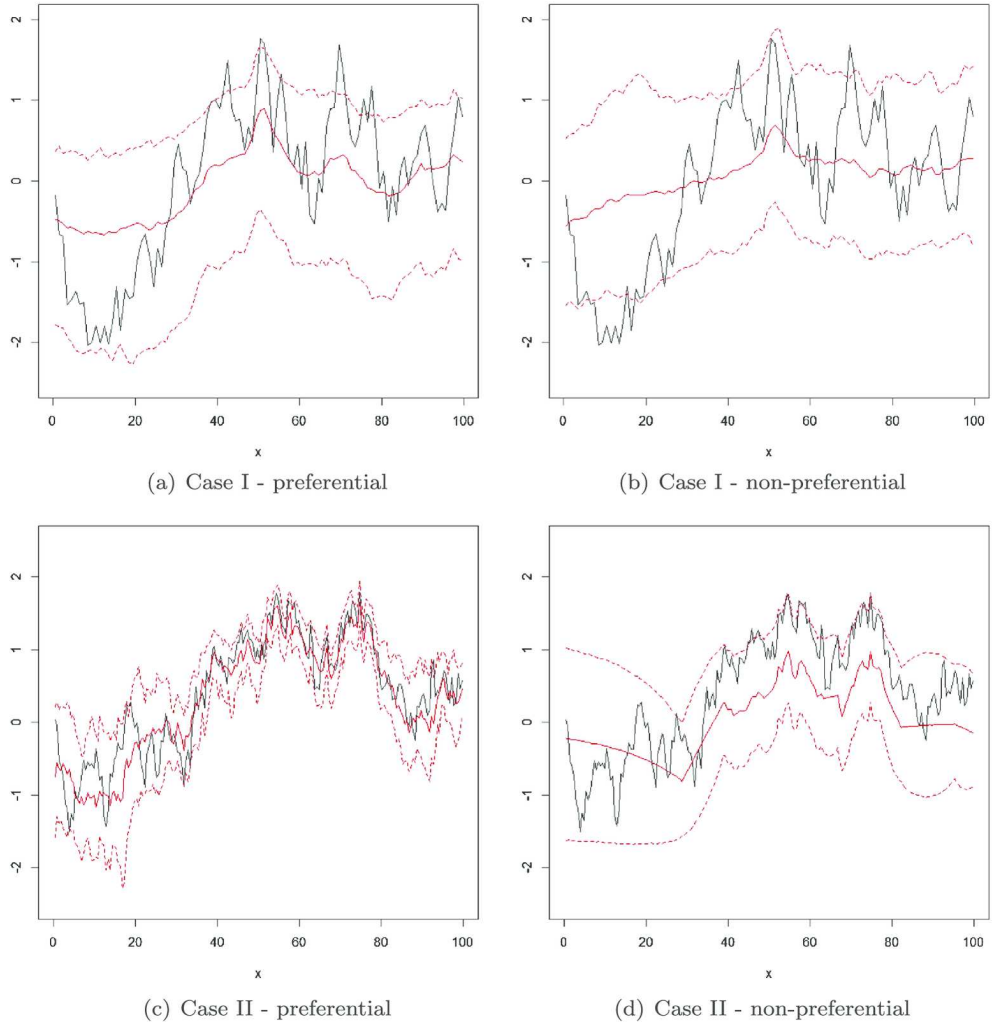}
\caption{Simulated process (solid line), posterior median
of the predictive distribution (red line) and the respective 95\%
credibility intervals (dashed lines) for $S$ obtained by considering
(left) and by not considering (right) the effect of preferential
sampling in {\it Cases I} and $II$.}\vspace*{-2pt}
\end{figure}

\begin{figure}[!htb]
\centering
\includegraphics{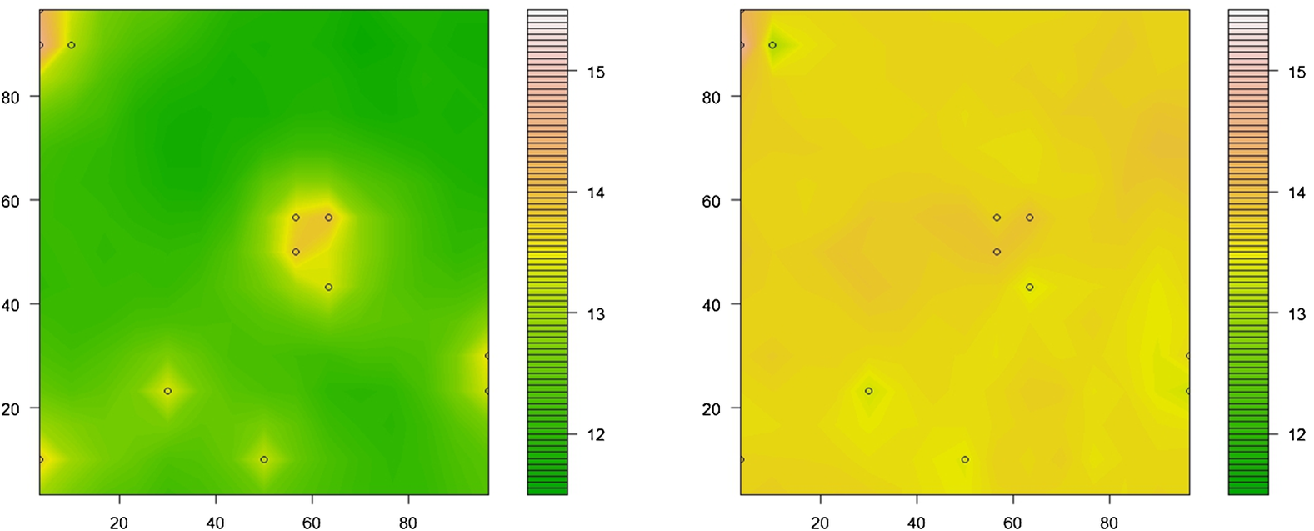}
\caption{Posterior mean of the predictive distribution of $[S+\mu]$
obtained considering (left) and without considering (right) the effect
of preferential sampling for {\it Case IV}.}
\end{figure}

All simulations presented in this section reflect situations where
there are few observed sample points. This is not uncommon in practice,
especially in monitoring studies of environmental or climatic phenomena
which are rare or difficult to detect. However, even under these
conditions, the use of models assuming the existence of preferential
sampling effects has produced variogram and kriging surfaces generally
closer to the true values when compared to the estimates produced
without considering this effect. Another major advantage of using these
models is the correction that is made during the inference about the
underlying process mean. Finally, we showed the ability\vadjust{\eject}
of this model to identify areas where the underlying process takes
extreme values, even in situations where there are no samples  nearby.

\citet{diggle2010geostatistical} evaluated the influence of
preferential sampling influence in the predictions of Gaussian
processes using $[S|{\bf y},\hat{\theta}]$ after correcting the bias
caused by preferential sampling. However, these predictions did not
take into account the information provided by {\bf x}. In the context
of environmental monitoring, \citet{shaddick2014case} also presented a
methodology for correcting the bias caused by a selective reduction of
sample sites. In contrast, the Bayesian approach provides samples
directly of the distribution of interest $[S|{\bf y},{\bf x}]$. The
comparison between the predictions assuming that $\theta$ is known
reinforced the conclusion that methods based on corrections of the
variogram's bias are not sufficient to reproduce the true uncertainty
associated with the predictive distribution of the underlying process
$S$.

Performing different simulated scenarios of prediction under
preferential sampling effects, \citet{gelfand2012effect} also concluded
that they affect more significantly the spatial prediction than the
estimation of the model parameters. In their paper, they discussed ways
to evaluate the effects of preferential sampling by comparing two
predicted surfaces. One of the forms of global comparison they
mentioned is associated with the local and the global mean squared
prediction error. The {\it Local Prediction Error} associated with
$x_0$, denoted $LPE(x_0)$, is given by
$$LPE(x_0)=E[\hat{S}(x_0)-S(x_0)]^2,$$
where $\hat{S}(x_0)$ is the predictor of $S$ in $x_0$. The {\it Global
Prediction Error} is given by
$$GPE=\frac{1}{|D|}\int_{D} LPE(x)dx.$$
Table 1 presents the $GPE$ values for each of the simulations. Based on
this table, one can observe a significant reduction obtained by
considering the effects of preferential sampling.

\begin{table}
\begin{tabular}{|c|c|c|}
\hline
Simulation & No preferential sampling & Preferential sampling\\
\hline \hline
{\it Case I} & 0.9496 & 0.7301 \\
\hline
{\it Case II} & 0.5533 & 0.1783 \\
\hline
{\it Case III} & 0.5286 & 0.4512 \\
\hline
{\it Case IV} & 2.1336 & 1.6789 \\
\hline
{\it Case V} & 1.6214 & 1.3482 \\
\hline
\end{tabular}
\caption{Global Prediction Error (GPE) for each of the simulations.}
\end{table}

The one-dimensional simulation showed that the $LPE$s are reduced when
they refer to the locations where the underlying process $S$ has lower
values. Similar conclusions were obtained in the second simulation,
since errors remain smaller in regions where the magnitude of $S$ is
lower when the preferential sampling effect is taken into account in
the modelling. Further simulations without the effect of preferential
sampling were also performed. In such cases, in general, the model
assuming preferential sampling produced posterior distributions of
$\beta$ centred at zero.

\subsection{Optimal Design}

Using the one-dimensional simulated data presented in {\it Case I}, 83
auxiliary points, required for evaluation of the utility function
(\ref{eq:utilidade1}) described in Section 3, were also used to form a
grid and obtaining the predictive variance reductions. Using the
samples of $\theta$, which were obtained from the posterior
distribution in the MCMC algorithm, we can generate samples of $d$ as
mentioned in Section 3 to obtain a new optimal sample point for the
case where $m=1$.

Figure 6 shows the histograms of the posterior pseudo-distribution of
$d$ with and without the preferential sampling assumption,
respectively. It can be noted that the optimal design choice without
preferential sampling leads the researcher to select locations where
there are no nearby points, i.e. in the interval $[5, 35]$. On the other
hand, under preferential sampling, the results lead the researcher to a
different direction. In this case, except in the subregions where there
are several observed samples, the other choices have similar expected
utilities. In summary, under preferential sampling, the optimal design
choice is notably changed when the researcher wants to reduce the
predictive variance.

\begin{figure}[!htb]
\centering
\includegraphics{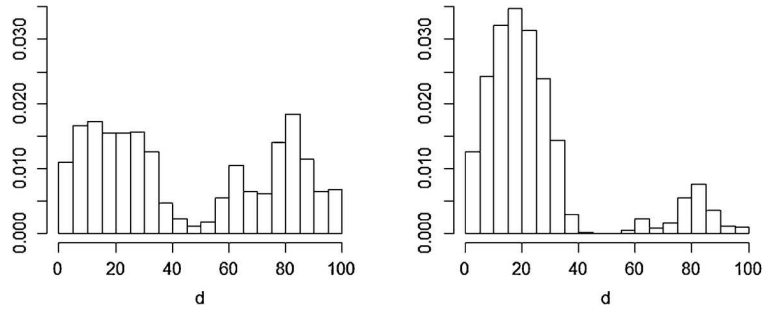}
\caption{Posterior pseudo-distributions of $d$ under the effect
of preferential sampling (left) and without considering this effect
(right) for {\it Case I}.}
\end{figure}

To illustrate the effect of the preferential sampling where $m>1$, we
obtained the pseudo-distribution of ${\bf d}=(d_1,d_2)'$ assuming that
two new locations would need to be added in the sample design. The
results are presented in Figure 7. Analysing the results for the
non-preferential model, it can be seen that the optimal solution could
be to add one point $d_1$ from interval $[0, 20]$ and to make no
restrictions to the location of the other point $d_2$, since $h({\bf
d})$ does not vary much. On the other hand, the results for the
preferential model indicate that the optimal design should include one
point $d_1$ from interval $[0, 4]$ and another point $d_2$ from interval
$[88, 92]$.

\begin{figure}[!htb]
\centering
\includegraphics{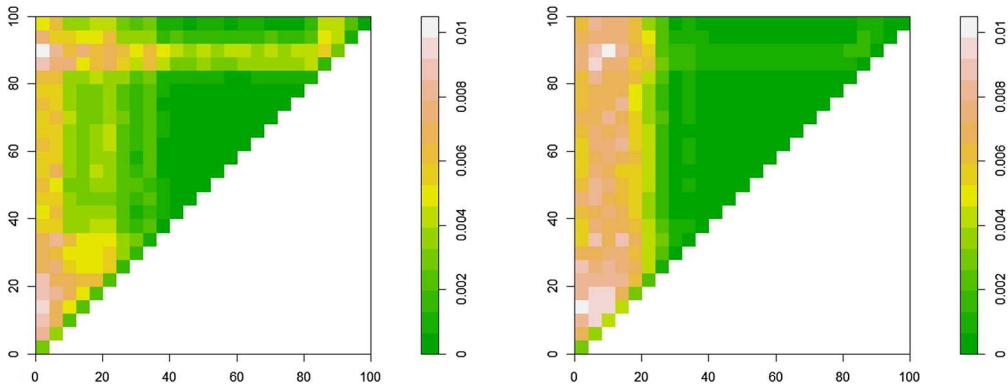}
\caption{Posterior pseudo-distributions of ${\bf d}=(d_1,d_2)'$ under
the effect of preferential sampling (left) and without considering this
effect (right) for {\it Case I}.}
\end{figure}

In the two-dimensional {\it Case IV}, 900 auxiliary points were also
used to form a grid and obtaining the predictive variance reductions,
as in the previous section. Figure 8 shows the posterior
pseudo-distributions of $d$ under preferential sampling and without
considering this effect for the case where $m=1$.

\begin{figure}[!htb]
\centering
\includegraphics{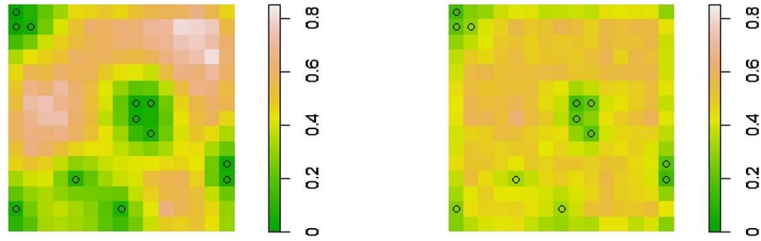}
\caption{Posterior pseudo-distribution of $d$ under the effect of preferential
sampling (left) and without considering this effect (right) for {\it
Case IV}. The densities of these pseudo-distribution are multiplied by
100 for better visualization.}
\end{figure}

It can be noted in Figure 8 that the areas with the highest expected
utility are not the same for the two models. As expected, the results
obtained without considering the effect of preferential sampling makes
the researcher choose locations far from the observed points. In
contrast, under preferential sampling, the largest utilities expected
are more dispersed over region $D$. Again, under preferential sampling,
the optimal design choice was changed, since the regions with few
sample sites also provide useful information about $S$.

The two simulated studies have involved situations where the inference
step has produced similar results (in the one-dimensional case) and
different results (in the two-dimensional case). Surprisingly, even in
the one-dimensional case, where the estimated variograms were similar,
the process of choosing the optimal design led to quite different
results.

However, the use of a different utility function can produce even more
extreme results. Even though the utility function chosen in this paper
was designed for our Geostatistics goals, other functions could be
considered. According to the results obtained from the simulated
studies, utility functions that depend on the underlying process $\mu$
may also be affected by the preferential sampling effect. This
situation will be explored in next section. There are other effects
that may affect the results, such as the choice of the auxiliary grid
(used to evaluate the predictive variance reduction) and the
discretization level of $D$. Current computational costs associated
with this methodology are still a barrier for a more thorough
evaluation of the degree of influence of each of these marginal
effects.

\subsubsection{Effectiveness of the Optimal Decision}

After obtaining the optimal design, one can evaluate if the results are
better under preferential sampling. Thus, we performed an analysis of
the $GPE$ after this decision to the two-dimensional simulated data
{\it Case IV}. The coordinates of the optimal design was
$x_d=(90.00;76.66)$, under preferential sampling, and
$x_d=(30.00;50.00)$ without this effect. In addiction, it was assumed
that $\tau^2=0$.

Finally, we proceed to the inference via classical Geostatistics
methods but including the optimal data point $Y_d$ obtained under
preferential sampling. The results provided a $GPE=1.8392$, which is
lower than that obtained using the optimal data point $Y_d$ pointed out
by the non-preferential model ($GPE=2.3348$). Thus, the optimal design
obtaining under preferential sampling was more advantageous even when
the inference is performed via classical Geostatistics methods.

\section{Case Study: Rainfall Data in Rio de Janeiro}

The methodology is now applied to a real scenario in the context of
monitoring networks. More specifically, we will analyse pluviometric
precipitation data obtained from 32 monitoring stations located in the
city of Rio de Janeiro, Brazil. The data refer to the period from 1 to
31 October 2005 and were obtained from the {\it Pereira Passos}
Institute, an official agency associated with the local government. The
rainfall during this month, which begins the rainy season in the
Brazil's Southeast, is of particular interest to meteorologists and
government agencies \citep{alves2005}. Figure 9 shows the map of the
Rio de Janeiro city with the respective precipitation levels observed.

\begin{figure}[!htb]
\centering
\includegraphics{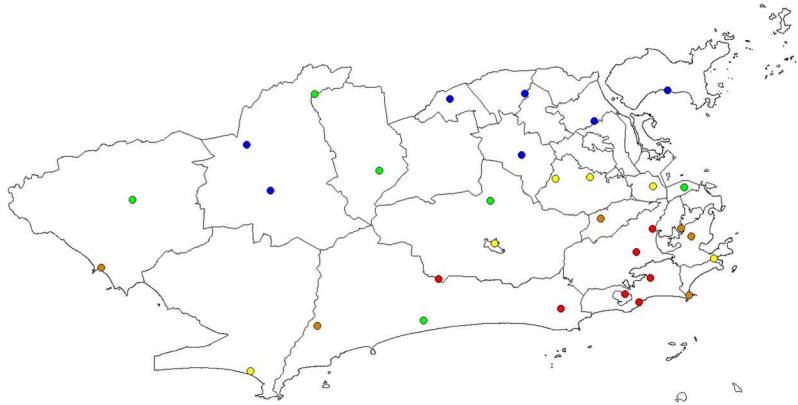}
\caption{Pluviometric precipitation in Rio de Janeiro city in October/2005
(separated according the 0.20, 0.40, 0.60, and 0.80 quantiles and
grouped by the colours: blue, green, yellow, brown, and red,
respectively).}
\end{figure}

Analysis of the spatial distribution of the precipitation seems to
indicate that the stations are more concentrated near places where
rainfall level is higher. Even though the geography and the spatial
distribution of economic activities could be considered as possible
causes of this sample design, the methodology for choosing a new design
point under the effect of preferential sampling can be employed here.

For inference, we used the same priors of the previous simulations (by
changing some hyperparameters) and the study area was partitioned into
$M=332$ subregions. We monitored 100,000 iterations in the MCMC
algorithm for both models and the first 10,000 were considered as
burn-in. The convergence of the chains was assessed by visual
inspection of several chains generated from different initial values.
Table 2 presents summaries of the posterior distributions for all model
parameters. An analysis of the results suggests that the effect of
preferential sampling is significant, indicating a positive association
between the sample design and the rainfall intensity.

\begin{table}[p]
\begin{tabular}{|c|c|c|}
\hline
Model parameters & Preferential & Non-Preferential\\
\hline \hline
$\tau^2$ & 1.25 (0.49; 2.90) & 1.32 (0.51; 3.45) \\
\hline
$\sigma^2$ & 4289.92 (2096.31; 10528.91) & 4132.72 (2120.30; 8514.38)\\
\hline
$\mu$ & 104.84 (97.73; 110.60) & 119.88 (111.49; 130.32) \\
\hline
$\phi$ & 10.69 (4.27; 26.75) & 10.43 (4.51; 22.76) \\
\hline
$\alpha$ & $-$3.84 ($-$4.24; $-$3.48) & ---\\
\hline
$\beta$ & 0.008 (0.002; 0.014) & ---\\
\hline
\end{tabular}
\caption{Posterior mean and 95\% credibility interval (in parentheses)
for model parameters.}
\end{table}

\begin{figure}[p]
\centering
\includegraphics{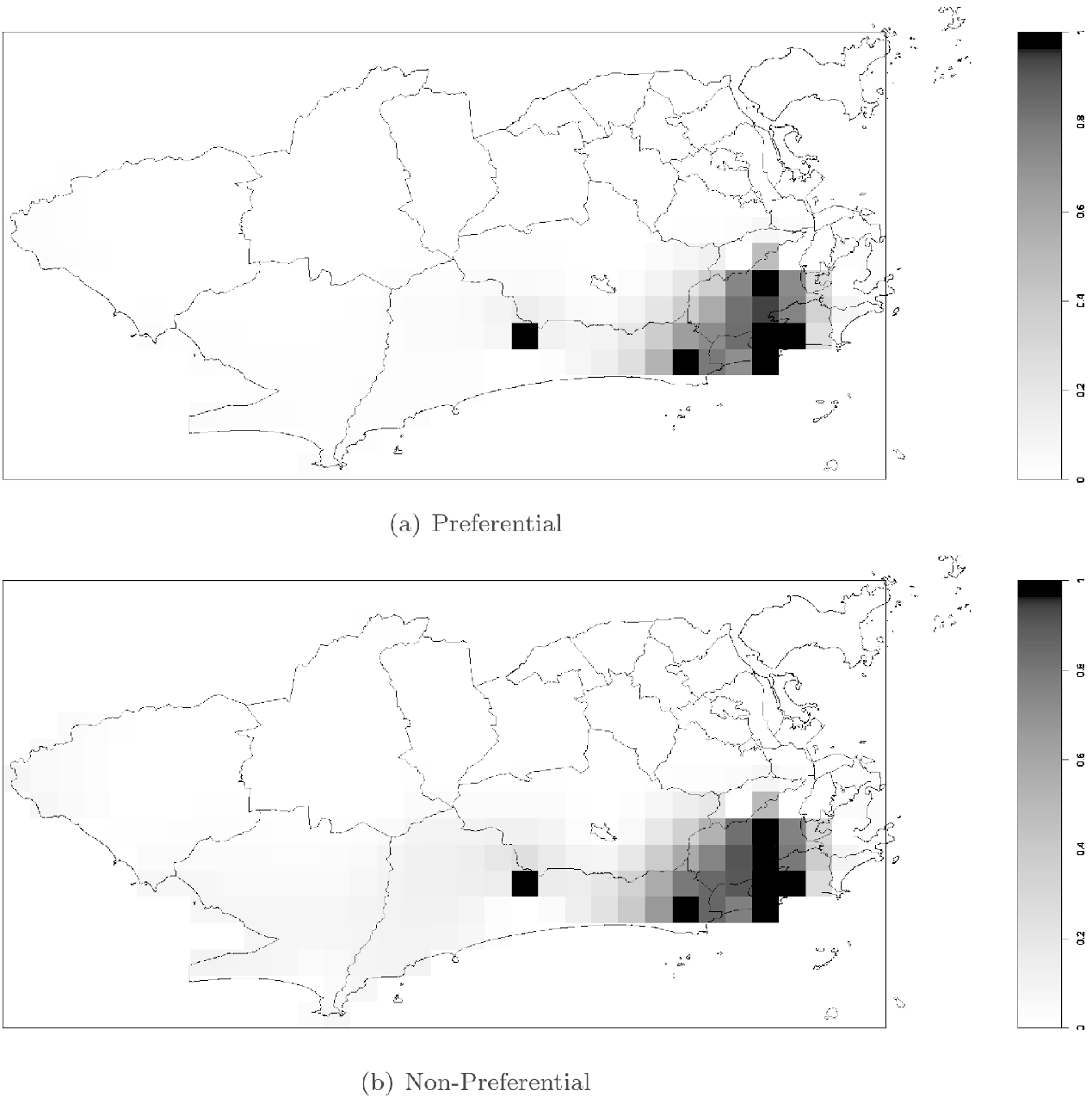}
\caption{Expected utilities $U(d)$ obtained under preferential sampling (a)
and without considering this effect (b) for the Rainfall data.}
\end{figure}

As for the other model parameters, we observed differences between the
estimates for the mean $\mu$, which can be explained by the presence of
preferential sampling effects. The utility function
(\ref{eq:utilidade2}) was used to obtaining the optimal design, with
$x_0 = 200$ mm. This utility function assigns greater utility for
regions where it is most likely to observe a monthly rainfall above 200
mm. The expected utility $U(d)$ obtained for each model is shown in
Figure 10.

The preferential model has concentrated high expected utility into a
small region in map since the utility function favours regions with the
highest probability of extremes values and due to the overestimation of
$\mu$. On the other hand, the non-preferential model has produced
expected utilities more spread out over the southern part of the
city.\looseness=-1

\section{Discussion}

The results obtained in Sections 5 and 6 help us to understand the
effects of preferential sampling on inference and optimal design choice
in the context of Geostatistics. The results of Section 5 showed that
the bias corrections made during the estimation step are not enough to
ensure a satisfactory prediction. On the other hand, the Bayesian
approach presented here allows us to make predictions about any
functional of $S$ directly from $[S\mid {\bf y,x}]$, which is the
correct predictive distribution under the effect of preferential
sampling. Furthermore, the inference about the parameters that define
the effect of preferential sampling, i.e. $\alpha$ and $\beta$, was
quite satisfactory in all simulations. Despite the high  associated
computational cost, the approach is computationally feasible and showed
satisfactory results.

In practical situations, it may be unlikely to assume that the design
is governed by a log-Gaussian Cox process. However, this model seems to
be flexible to obtain understanding about the consequences if the
researcher has no covariates or a better explanation for its true
causes. Although not applicable in the strictly spatial context, an
alternative approach to detect and deal with preferential sampling is
spatio-temporal analysis, since the changes on site locations over time
may be informative about the association between its configuration and
the underlying process. \citet{zidek2014reducing} present a method that
can learn about the preferential selection process over time and that
allows the researcher to deal with its effects.

Traditional approaches to deal with preferential sampling also include
the fit of trend surfaces to assess first order effects and, in the
context of survey sampling methods, weighting schemes. However, both
require some prior knowledge about the sample selection processes. The
researcher must carefully evaluate the (theoretical or empirical)
evidence of preferential sampling effects. The use of this approach can
produce misleading results otherwise. If there is evidence of
preferential sampling, the authors believe the use of models based on
latent processes, such as log-Gaussian Cox process, should be used even
when the dependence structure between $X$ and $S$ is not completely
known.

As evidenced by simulations, the effects of preferential sampling on
optimal design choice cannot be disregarded. The knowledge about the
sampling pattern reduces the predictive variance of $S$ in areas poorly
sampled, substantially changing the optimal decision. On the other
hand, a utility function based on exceedances seems to overestimate
$\mu$ under preferential sampling. As a result, the ideal region to
receive a new sample becomes very small. However, obtaining the optimal
design in situations where there is a need to monitor extreme events or
exceedances may be not simple. \citet{chang2007designing} presented an
approach to deal with some challenges arising in designing networks for
monitoring fields of extremes, as the loss of spatial dependence and
the limitations of conventional approaches.

Potential areas to apply this methodology include the study of
phenomena scarcely observed and those in which, due to researcher's
interest or limited resources, can only be observed in locations
considered critical. One can include the monitoring of mosquitoes or
other disease spreading pests that only become detectable in places
where its occurrence is very high among the phenomena that are scarcely
observed and difficult to be detected. Another potential application is
the monitoring of maximum and minimal temperatures in regions near
airports or industrial plants. In both cases, it seems reasonable to
infer that the way the phenomenon is observed can be related to the
underlying process. In this situation, the optimal design choice can
change significantly, since the spatial point pattern brings valuable
information to the researcher.

The methodology employed can be adjusted in order to produce
pseudo-distributions of $d$ more peaked around the mode. The strategy
based on simulated annealing \citep{muller1999simulation} can be
explored for this purpose. The approach of \citet{muller1999simulation}
for optimal design choice has computational advantages in comparison
with traditional procedures of optimization. The choice of a utility
function based on predictive variance reduction is easily justified in
Geostatistics. Under preferential sampling, the use of utility
functions that directly depend on $\mu$ seems to be more affected than
those based in variance reductions.\looseness=-1

This methodology also has an expensive computational cost and
alternative procedures can be used to deal with this problem, e.g.
\citet{simpson2011going} that use the {\it Integrated Nested Laplace
Aproximation} -- INLA methods \citep{rue2009approximate}, the use of
Predictive Process \citep{banerjee2008gaussian}, methods to approximate
likelihoods
\citep[see][]{stein2004approximating,fuentes2007approximate} and the
use of sparse covariance matrices
\citep{furrer2006covariance}.\looseness=-1

The authors recognize the limitations in the simulation studies
presented here. The complexity and the variety of situations arising
from designing under preferential sampling lead us to focus on specific
features rather than obtaining more general conclusions. Finally, the
authors also invite researchers interested in reproducing the
methodology presented in this paper to contact the authors in order to
obtain more details about the computational implementation.

\appendix
\section*{Appendix A}

In the preferential model, the full conditional distributions of
$S,\beta$ and $\alpha$ are proportional to
\begin{align*}
p(S \mid \mu,\tau^{-2},\sigma^{-2},\phi,\alpha,\beta,{\bf x,y})
& \propto\exp \left\{ -\frac{1}{2\tau^2} [S_{\bf y}'S_{\bf y} - 2S_{\bf y}'({\bf y}
- \mu {\bf 1}) ]  + \beta S'{\bf n} - \frac{S'R_{M}^{-1}S}{2\sigma^2} \right\}\\
&\quad \times \exp \left\{ - \Delta e^{\alpha} \sum^{M} \exp(\beta S(x_i)) \right\},
\end{align*}
\begin{align*}
p(\beta \mid S, \mu,\tau^{-2},\sigma^{-2},\phi,\alpha,{\bf x,y}) &\propto
p({\bf x} \mid S,\alpha,\beta)p(\beta)\\
&\propto \exp \left\{ \beta S'{\bf n} - \Delta e^{\alpha} \sum^{M}
\exp(\beta S(x_i)) - \frac{\beta^2}{2k} \right\} ,
\end{align*}
\[
p(\alpha \mid S, \mu,\tau^{-2},\sigma^{-2},\phi,\beta,{\bf x,y}) \propto \exp
\left\{ n\alpha - \Delta e^{\alpha} \sum^{M} \exp(\beta S(x_i))
- \frac{\alpha^2}{2k}  \right\}.
\]
In the MCMC, these quantities can be updated in Metropolis
steps assuming a Gaussian proposal distribution centred in the previous
values sampled. Then, we accept the proposal with probability
\begin{align*}
p_S&=\exp \biggl\{ -\frac{1}{2\tau^2} [S_{\bf y}'^{prop}S_{\bf y}^{prop} -
S_{\bf y}'S_{\bf y} - 2(S_{\bf y}^{prop}-S_{\bf y})'({\bf y} - \mu {\bf 1}) ]  +
\beta (S^{prop}-S)'{\bf n} \biggr\} \\
&\quad \times \exp \biggl\{ \Delta e^{\alpha} \sum^{M} [\exp(\beta S(x_i))
-\exp(\beta S(x_i)^{prop})] + \frac{(S'R_{M}^{-1} S
- S'^{prop} R_{M}^{-1} S^{prop})}{2\sigma^2} \biggr\},
\end{align*}
\[
p_{\beta}=\exp \biggl\{ (\beta^{prop}\,{-}\,\beta) S'{\bf n} \,{+}\,
\Delta e^{\alpha} \sum^{M} (\exp(\beta S(x_i))\,{-}\,\exp(\beta^{prop} S(x_i))) \,{+}\,
\frac{(\beta^2-\beta^{2prop})}{2k} \biggr\},
\]
\[
p_{\alpha}=\exp \biggl\{ (\alpha^{prop}-\alpha) n + (e^{\alpha}-e^{\alpha^{prop}})
\Delta \sum^{M} \exp (\beta S(x_i)) + \frac{(\alpha^2-\alpha^{2prop})}{2k} \biggr\},
\]
respectively. The vector ${\bf n}'=(n_1,n_2,\dots,n_M)$
represents the number of observations in each subregion, where
$\sum_{i=1}^M n_i=n$.

In both models, the full conditional distribution of $\phi$ is
proportional to
\[
p(\phi \mid S,\mu,\tau^{-2},\sigma^{-2},\alpha,\beta,{\bf
y}) \propto |R_M|^{-1/2} \phi^{a_{\phi}-1} \exp \left\{
-\frac{S'R_{M}^{-1}S}{2\sigma^2} -b_{\phi}\phi  \right\}
\]
and this parameter can be updated in Metropolis steps
assuming the following proposal distribution
$$q(\phi^{prop} \mid \phi) \sim Lognormal\left( \ln(\phi) - \delta/2; \delta \right),$$
where $\delta$ must be chosen in order to produce reasonable
acceptance rates in MCMC. Then, we accept the proposal with probability
\begin{align*}
p_\phi&= \left( \frac{|R_M|^{prop}}{|R_M|} \right)^{-1/2}
\left( \frac{\phi^{prop}}{\phi}\right)^{a_{\phi}}
\times \exp \biggl\{ -\frac{(S'R_{M}^{prop-1}S-S'R_{M}^{-1}S)}{2\sigma^2}
+b_{\phi}(\phi-\phi^{prop})\\
&\quad - \frac{(\ln\phi-\ln\phi^{prop}+ \delta/2)^2-(\ln\phi^{prop}
-\ln\phi+ \delta/2)^2}{2\delta}\biggr\}.
\end{align*}

\section*{Appendix B}

To evaluate the integral in expression of $u({\bf d}, \theta, {\bf
y_d})$ a discretization of the region $D$ can be applied yielding $M$
subregions, as described in Section 2. Then, we have that $u({\bf d},
\theta, {\bf y_d})$ can be approximated by
$$\tilde{u}({\bf d},\theta,{\bf y_d}) = \frac{1}{M} \sum_{i}
[V(S_i \mid \theta,{\bf y})-V(S_i \mid \theta, {\bf y}, {\bf y_d})],$$
where $[S_i \mid \theta,{\bf y}]$ and $[S_i \mid \theta, {\bf
y}, {\bf y_d}]$, $i =1,\dots,M$, are Gaussian with respective means
\begin{align*}
\sigma^2{\bf r}^{'}_{n}(\tau^2 I_{n} + \sigma^2R_{n})^{-1}({\bf
y}-{\bf 1}\mu) \quad \mbox{and} \quad \sigma^2{\bf r}^{'}_{n+m}(\tau^2
I_{n+m} + \sigma^2R_{n+m})^{-1}({\bf y}^{*}-{\bf 1}\mu)
\end{align*}
and variances
\begin{align*}
\sigma^2-\sigma^2{\bf r}^{'}_{n}(\tau^2 I_{n} +
\sigma^2R_{n})^{-1}\sigma^2{\bf r}_{n} \quad \mbox{and} \quad
\sigma^2-\sigma^2{\bf r}^{'}_{n+m}(\tau^2 I_{n+m} +
\sigma^2R_{n+m})^{-1}\sigma^2{\bf r}_{n+m},
\end{align*}
where ${\bf y}^{*}=({\bf y}, {\bf y_d})'$.

\section*{Appendix C}

Following the procedure proposed in Section 4, we can make an
approximation in $p({\bf x} \mid S, \alpha, \beta)$ to obtain an
analytical expression of $V[S \mid \theta,{\bf y},{\bf x}]$. In
particular, expanding the exponential function in Taylor series around
zero, up to the second order term, we obtain the following
approximation
$$\sum_i \exp(\alpha+\beta S(x_i)) \approx e^{\alpha}\left({\bf 1_{M}}
+  \beta{\bf 1_{M}}'S + \frac{\beta^2}{2}S'S \right).$$

By inserting this expression in $p({\bf x} \mid S, \alpha, \beta)$, it
can be shown that the conditional distribution $[S \mid \theta,{\bf
y},{\bf x}]$ becomes Gaussian with mean vector $\Theta$ and covariance
matrix $\Sigma$ given by
\begin{align*}
\Theta=\Sigma \times \left(\begin{array}{c} ({\bf y_n}-\mu {\bf
n})\tau^{-2} + \beta {\bf n} - \Delta \beta e^{\alpha}{\bf n} \\[3pt] -
\Delta \beta e^{\alpha}{\bf 1_N} \end{array} \right)\quad  \mbox{and}
\end{align*}
$$\Sigma= \left(\begin{array}{cc}
(\tau^{-2} + \Delta \beta^2 e^{\alpha})I_{{\bf n}} + R_{n}^{-1}R_{n,N} A^{-1}
R_{N,n}R_{n}^{-1} + \sigma^{-2}R_{n}^{-1} & -R_{n}^{-1}R_{n,N} A^{-1} \\[3pt]
-A^{-1}R_{N,n}R_{n}^{-1}  & \Delta \beta^{2} e^{\alpha} I_{N} + A^{-1} \end{array}\right)^{-1},$$
where $A = \sigma^2 R_N- \sigma^2 R_{N,n} R_{n}^{-1}
R_{n,N}$, and the vectors {\bf n} and ${\bf y_n}$ represent the number
of observations and the total observed in each subregion of $D$,
respectively. If \mbox{$\beta=0$}, this matrix becomes equal to the kriging
variance matrix traditionally obtained in Geostatistics.

\begin{acknowledgement}
  This paper is based on the doctoral thesis of the first author under the
supervision of the second author. The authors thank the referees for
numerous comments and suggestions that led to substantial improvement
of the presentation. The second author wishes to thank the support from
CNPq-BRAZIL.
\end{acknowledgement}


\begin{thebibliography}{38}

\bibitem[{Alves et~al.(2005)Alves, Marengo, J\'{u}nior, and Castro}]{alves2005}
Alves, L., Marengo, J., J\'{u}nior, H., and Castro, C. (2005).
\newblock \enquote{Beginning of the rainy season in southeastern Brazil: Part 1
 -- Observational studies (in Portuguese).}
\newblock \emph{Revista Brasileira de Meteorologia\/}, 20(3): 385--394.
\bptok{addids}%
\endbibitem

\bibitem[{Banerjee et~al.(2008)Banerjee, Gelfand, Finley, and
 Sang}]{banerjee2008gaussian}
Banerjee, S., Gelfand, A.~E., Finley, A.~O., and Sang, H. (2008).
\newblock \enquote{Gaussian predictive process models for large spatial data
 sets.}
\newblock \emph{Journal of the Royal Statistical Society: Series B (Statistical
 Methodology)\/}, 70(4): 825--848.
\bid{doi={10.1111/j.1467-9868.2008.00663.x}, issn={1369-7412},
mr={2523906}}
\bptok{addids}%
\endbibitem

\bibitem[{Bornn et~al.(2012)Bornn, Shaddick, and Zidek}]{bornn2012}
Bornn, L., Shaddick, G., and Zidek, J. (2012).
\newblock \enquote{Modeling nonstationary processes through dimension
 expansion.}
\newblock \emph{Journal of the American Statistical Association\/}, 107(497):
 281--289.
\bid{doi={10.1080/01621459.\\2011.646919}, issn={0162-1459},
mr={2949359}}
\bptok{addids}%
\endbibitem

\bibitem[{Boukouvalas et~al.(2009)Boukouvalas, Cornford, and
 Stehlik}]{boukouvalas2009}
Boukouvalas, A., Cornford, D., and Stehlik, M. (2009).
\newblock \enquote{Approximately optimal experimental design for
 heteroscedastic Gaussian process models.}
\newblock Technical report, Nonlinear Complexity Research Group (NCRG).
\bptok{addids}%
\endbibitem

\bibitem[{Chang et~al.(2007)Chang, Fu, Le, and Zidek}]{chang2007designing}
Chang, H., Fu, A.~Q., Le, N.~D., and Zidek, J.~V. (2007).
\newblock \enquote{Designing environmental monitoring networks to measure
 extremes.}
\newblock \emph{Environmental and Ecological Statistics\/},\break 14(3): 301--321.
\bid{doi={10.1007/s10651-007-0020-5}, issn={1352-8505}, mr={2405332}}
\bptok{addids}%
\endbibitem

\bibitem[{Cressie(1993)}]{cressie1993statistics}
Cressie, N. (1993).
\newblock \emph{Statistics for spatial data\/}.
\newblock John Wiley and Sons, Inc.
\bid{mr={1239641}}
\bptok{addids}%
\endbibitem

\bibitem[{DeGroot(2005)}]{degroot2005optimal}
DeGroot, M.~H. (2005).
\newblock \emph{Optimal statistical decisions\/}, volume~82.
\newblock Wiley-Interscience.
\bid{doi={10.1002/0471729000}, mr={2288194}}
\bptok{addids}%
\endbibitem

\bibitem[{Diggle(1983)}]{diggle1983}
Diggle, P. (1983).
\newblock \emph{Statistical Analysis of Spatial Point Patterns\/}.
\newblock London: Academic Press.
\bid{mr={0743593}}
\bptok{addids}%
\endbibitem

\bibitem[{Diggle and Lophaven(2006)}]{diggle2006bayesian}
Diggle, P.~J. and Lophaven, S. (2006).
\newblock \enquote{Bayesian geostatistical design.}
\newblock \emph{Scandinavian Journal of Statistics\/}, 33(1): 53--64.
\bid{doi={10.1111/\\j.1467-9469.2005.00469.x}, issn={0303-6898},
mr={2255109}}
\bptok{addids}%
\endbibitem

\bibitem[{Diggle et~al.(2010)Diggle, Menezes, and
 Su}]{diggle2010geostatistical}
Diggle, P.~J., Menezes, R., and Su, T.-L. (2010).
\newblock \enquote{Geostatistical inference under preferential sampling.}
\newblock \emph{Journal of the Royal Statistical Society: Series~C (Applied
 Statistics)\/}, 59(2): 191--232.
\bid{doi={10.1111/j.1467-9876.\\2009.00701.x}, issn={0035-9254},
mr={2744471}}
\bptok{addids}%
\endbibitem

\bibitem[{Diggle and Ribeiro(2007)}]{diggle2007model}
Diggle, P.~J. and Ribeiro, P.~J. (2007).
\newblock \emph{Model-based geostatistics\/}.
\newblock Springer.
\bid{mr={2293378}}
\bptok{addids}%
\endbibitem

\bibitem[{Diggle et~al.(1998)Diggle, Tawn, and Moyeed}]{diggle1998model}
Diggle, P.~J., Tawn, J., and Moyeed, R. (1998).
\newblock \enquote{Model-based geostatistics.}
\newblock \emph{Journal of the Royal Statistical Society: Series C (Applied
 Statistics)\/}, 47(3): 299--350.
\bid{doi={10.1111/1467-9876.00113}, issn={0035-9254}, mr={1626544}}
\bptok{addids}%
\endbibitem

\bibitem[{Ding et~al.(2008)Ding, Rosner, and M{\"u}ller}]{ding2008bayesian}
Ding, M., Rosner, G.~L., and M{\"u}ller, P. (2008).
\newblock \enquote{Bayesian optimal design for phase II screening trials.}
\newblock \emph{Biometrics\/}, 64(3): 886--894.
\bid{doi={\\10.1111/j.1541-0420.2007.00951.x}, issn={0006-341X},
mr={2526640}}
\bptok{addids}%
\endbibitem

\bibitem[{Fern{\'a}ndez et~al.(2005)Fern{\'a}ndez, Real, Couto, Aboal, and
 Carballeira}]{fernandez2005effect}
Fern{\'a}ndez, J., Real, C., Couto, J., Aboal, J., and Carballeira, A.
(2005).
\newblock \enquote{The effect of sampling design on extensive bryomonitoring
 surveys of air pollution.}
\newblock \emph{Science of the total environment\/}, 337(1): 11--21.
\bptok{addids}%
\endbibitem

\bibitem[{Fuentes(2002)}]{fuentes2002spectral}
Fuentes, M. (2002).
\newblock \enquote{Spectral methods for nonstationary spatial processes.}
\newblock \emph{Biometrika\/}, 89(1): 197--210.
\bid{doi={10.1093/biomet/89.1.197}, issn={0006-3444}, mr={1888368}}
\bptok{addids}%
\endbibitem

\bibitem[{Fuentes(2007)}]{fuentes2007approximate}
--- (2007).
\newblock \enquote{Approximate likelihood for large irregularly spaced spatial
 data.}
\newblock \emph{Journal of the American Statistical Association\/}, 102(477):
 321--331.
\bid{doi={10.1198/016214506000000852}, issn={0162-1459}, mr={2345545}}
\bptok{addids}%
\endbibitem

\bibitem[{Fuentes and Smith(2001)}]{fuentes2001new}
Fuentes, M. and Smith, R.~L. (2001).
\newblock \enquote{A new class of nonstationary spatial models.}
\newblock Technical report, North Carolina State University, Raleigh, NC.
\bptok{addids}%
\endbibitem

\bibitem[{Furrer et~al.(2006)Furrer, Genton, and Nychka}]{furrer2006covariance}
Furrer, R., Genton, M.~G., and Nychka, D. (2006).
\newblock \enquote{Covariance tapering for interpolation of large spatial
 datasets.}
\newblock \emph{Journal of Computational and Graphical Statistics\/}, 15(3).
\bid{doi={10.1198/106186006X132178}, issn={1061-8600}, mr={2291261}}
\bptok{addids}%
\endbibitem

\bibitem[{Gelfand et~al.(2012)Gelfand, Sahu, and Holland}]{gelfand2012effect}
Gelfand, A.~E., Sahu, S.~K., and Holland, D.~M. (2012).
\newblock \enquote{On the effect of preferential sampling in spatial
 prediction.}
\newblock \emph{Environmetrics\/}, 23(7): 565--578.
\bid{doi={10.1002/env.2169}, issn={1180-4009}, mr={3020075}}
\bptok{addids}%
\endbibitem

\bibitem[{Gumprecht et~al.(2009)Gumprecht, M{\"u}ller, and
 Rodr{\'\i}guez-D{\'\i}az}]{gumprecht2009designs}
Gumprecht, D., M{\"u}ller, W.~G., and Rodr{\'\i}guez-D{\'\i}az, J.~M.
(2009).
\newblock \enquote{Designs for detecting spatial dependence.}
\newblock \emph{Geographical Analysis\/}, 41(2): 127--143.
\bptok{addids}%
\endbibitem

\bibitem[{Higdon et~al.(1999)Higdon, Swall, and Kern}]{higdon1999non}
Higdon, D., Swall, J., and Kern, J. (1999).
\newblock \enquote{Non-stationary spatial modeling.}
\newblock In: \emph{Bayesian statistics\/} (J. M. Bernardo, J. O. Berger, A. P. Dawid, and A. F. M. Smith), 6(1): 761--768.
\bptok{addids}%
\endbibitem

\bibitem[{M{\o}ller et~al.(1998)M{\o}ller, Syversveen, and
 Waagepetersen}]{moller1998log}
M{\o}ller, J., Syversveen, A.~R., and Waagepetersen, R.~P. (1998).
\newblock \enquote{Log-Gaussian Cox processes.}
\newblock \emph{Scandinavian Journal of Statistics\/}, 25(3): 451--482.
\bid{doi={10.1111/1467-9469.00115}, issn={0303-6898}, mr={1650019}}
\bptok{addids}%
\endbibitem

\bibitem[{M{\o}ller and Waagepetersen(2007)}]{moller2007modern}
M{\o}ller, J. and Waagepetersen, R.~P. (2007).
\newblock \enquote{Modern statistics for spatial point processes.}
\newblock \emph{Scandinavian Journal of Statistics\/}, 34(4): 643--684.
\bid{doi={10.1111/j.1467-9469.2007.00569.x}, issn={0303-6898},
mr={2392447}}
\bptok{addids}%
\endbibitem

\bibitem[{M{\"u}ller(1999)}]{muller1999simulation}
M{\"u}ller, P. (1999).
\newblock \enquote{Simulation-based optimal design.}
\newblock In: \emph{Bayesian statistics\/} (J.~M.~Bernardo, J. O. Berger, A. P. Dawid, and A. F. M. Smith), 6: 459--474.
\bid{mr={1723509}}
\bptok{addids}%
\endbibitem

\bibitem[{M{\"u}ller et~al.(2007)M{\"u}ller, Berry, Grieve, Smith, and
 Krams}]{muller2007simulation}
M{\"u}ller, P., Berry, D.~A., Grieve, A.~P., Smith, M., and Krams, M.
(2007).
\newblock \enquote{Simulation-based sequential Bayesian design.}
\newblock \emph{Journal of Statistical Planning and Inference\/}, 137(10):
 3140--3150.
\bid{doi={10.1016/j.jspi.2006.\\05.021}, issn={0378-3758},
mr={2364157}}
\bptok{addids}%
\endbibitem

\bibitem[{M{\"u}ller et~al.(2004)M{\"u}ller, Sans{\'o}, and
 De~Iorio}]{muller2004optimal}
M{\"u}ller, P., Sans{\'o}, B., and De~Iorio, M. (2004).
\newblock \enquote{Optimal Bayesian design by inhomogeneous Markov chain
 simulation.}
\newblock \emph{Journal of the American Statistical Association\/}, 99(467):
 788--798.
\bid{doi={10.1198/\\016214504000001123}, issn={0162-1459},
mr={2090911}}
\bptok{addids}%
\endbibitem

\bibitem[{M{\"u}ller and Stehl{\'\i}k(2010)}]{muller2010compound}
M{\"u}ller, W.~G. and Stehl{\'\i}k, M. (2010).
\newblock \enquote{Compound optimal spatial designs.}
\newblock \emph{Environmetrics\/}, 21(3-4): 354--364.
\bid{doi={10.1002/env.\\1009}, issn={1180-4009}, mr={2842248}}
\bptok{addids}%
\endbibitem

\bibitem[{Pati et~al.(2011)Pati, Reich, and Dunson}]{pati2011bayesian}
Pati, D., Reich, B.~J., and Dunson, D.~B. (2011).
\newblock \enquote{Bayesian geostatistical modelling with informative sampling
 locations.}
\newblock \emph{Biometrika\/}, 98(1): 35--48.
\bid{doi={10.1093/biomet/asq067}, issn={0006-3444}, mr={2804208}}
\bptok{addids}%
\endbibitem

\bibitem[{Ripley(2005)}]{ripley2005spatial}
Ripley, B.~D. (2005).
\newblock \emph{Spatial statistics\/}, volume 575.
\newblock Wiley.
\bid{mr={0624436}}
\bptok{addids}%
\endbibitem

\bibitem[{Rue et~al.(2009)Rue, Martino, and Chopin}]{rue2009approximate}
Rue, H., Martino, S., and Chopin, N. (2009).
\newblock \enquote{Approximate Bayesian inference for latent Gaussian models by
 using integrated nested Laplace approximations.}
\newblock \emph{Journal of the Royal Statistical Society: Series B (Statistical
 Methodology)\/}, 71(2): 319--392.
\bid{doi={10.1111/j.1467-9868.2008.00700.x}, issn={1369-7412}, mr={2649602}}
\bptok{addids}%
\endbibitem

\bibitem[{Shaddick and Zidek(2014)}]{shaddick2014case}
Shaddick, G. and Zidek, J.~V. (2014).
\newblock \enquote{A case study in preferential sampling: Long term monitoring
 of air pollution in the UK.}
\newblock \emph{Spatial Statistics\/}, 9: 51--65.
\bptok{addids}%
\endbibitem

\bibitem[{Simpson et~al.(2011)Simpson, Illian, Lindgren, S{\o}rbye, and
 Rue}]{simpson2011going}
Simpson, D., Illian, J., Lindgren, F., S{\o}rbye, S., and Rue, H.
(2011).
\newblock \enquote{Going off grid: Computationally efficient inference for
 log-Gaussian Cox processes.}
\newblock arXiv:\href{http://arxiv.org/abs/1111.0641}{1111.0641}.
\bptok{addids}%
\endbibitem

\bibitem[{Stein et~al.(2004)Stein, Chi, and Welty}]{stein2004approximating}
Stein, M.~L., Chi, Z., and Welty, L.~J. (2004).
\newblock \enquote{Approximating likelihoods for large spatial data sets.}
\newblock \emph{Journal of the Royal Statistical Society: Series B (Statistical
 Methodology)\/}, 66(2): 275--296.
\bid{doi={10.1046/j.1369-\\7412.2003.05512.x}, issn={1369-7412},
mr={2062376}}
\bptok{addids}%
\endbibitem

\bibitem[{Stroud et~al.(2001)Stroud, M{\"u}ller, and
 Rosner}]{stroud2001optimal}
Stroud, J.~R., M{\"u}ller, P., and Rosner, G.~L. (2001).
\newblock \enquote{Optimal sampling times in population pharmacokinetic
 studies.}
\newblock \emph{Journal of the Royal Statistical Society: Series~C (Applied
 Statistics)\/}, 50(3): 345--359.
\bid{doi={10.\\1111/1467-9876.00239}, issn={0035-9254}, mr={1856330}}
\bptok{addids}%
\endbibitem

\bibitem[{Waagepetersen(2004)}]{waagepetersen2004convergence}
Waagepetersen, R. (2004).
\newblock \enquote{Convergence of posteriors for discretized log Gaussian Cox
 processes.}
\newblock \emph{Statistics \& Probability Letters\/}, 66(3): 229--235.
\bid{doi={10.1016/j.spl.2003.10.014}, issn={0167-7152}, mr={2044908}}
\bptok{addids}%
\endbibitem

\bibitem[{Zhu and Stein(2005)}]{zhu2005spatial}
Zhu, Z. and Stein, M.~L. (2005).
\newblock \enquote{Spatial sampling design for parameter estimation of the
 covariance function.}
\newblock \emph{Journal of Statistical Planning and Inference\/}, 134(2):
 583--603.
\bid{doi={10.1016/j.jspi.2004.04.017}, issn={0378-3758}, mr={2200074}}
\bptok{addids}%
\endbibitem

\bibitem[{Zidek et~al.(2014)Zidek, Shaddick, Taylor et~al.}]{zidek2014reducing}
Zidek, J.~V., Shaddick, G., Taylor, C.~G., et~al. (2014).
\newblock \enquote{Reducing estimation bias in adaptively changing monitoring
 networks with preferential site selection.}
\newblock \emph{The Annals of Applied Statistics\/}, 8(3): 1640--1670.
\bid{doi={10.\\1214/14-AOAS745}, issn={1932-6157}, mr={3271347}}
\bptok{addids}%
\endbibitem

\bibitem[{Zidek et~al.(2000)Zidek, Sun, and Le}]{zidek2000designing}
Zidek, J.~V., Sun, W., and Le, N.~D. (2000).
\newblock \enquote{Designing and integrating composite networks for monitoring
 multivariate Gaussian pollution fields.}
\newblock \emph{Journal of the Royal Statistical Society: Series C (Applied
 Statistics)\/}, 49(1): 63--79.
\bid{doi={10.1111/1467-9876.00179}, issn={0035-9254}, mr={1817875}}
\bptok{addids}%
\endbibitem

\end{thebibliography}
\end{document}


\begin{frontmatter}
\title{Rejoinder\thanksref{T1}}
%
%
%
\relateddois{T1}{Main article DOI: \relateddoi[ms=BA944]{Related
item:}{10.1214/15-BA944}.}
\runtitle{Rejoinder}

\begin{aug}
\author[a]{\fnms{Gustavo da Silva} \snm{Ferreira}\corref{}\ead
[label=e1]{gustavo.ferreira@ibge.gov.br}}
\and
\author[b]{\fnms{Dani} \snm{Gamerman}\ead[label=e2]{dani@im.ufrj.br}}


\runauthor{G. S. Ferreira and D. Gamerman}


\address[a]{National School of Statistical Sciences,
Brazilian Institute of Geography and Statistics,
Rio de Janeiro, Brazil, \printead{e1}}
\address[b]{Department of Statistical Methods,
Mathematical Institute, Federal University of Rio de Janeiro,
Rio de Janeiro, Brazil, \printead{e2}}

\end{aug}

%


\end{frontmatter}


\section{Introduction}

We would like to start by thanking the Editor of BA for the opportunity
to discuss our work and the discussants Peter J. Diggle, Michael G.
Chipeta, James V. Zidek, Noel Cressie and Raymond L. Chambers for a
very thorough evaluation of our contribution and for their thoughts on
the topic.
Our rejoinder to the discussion will be presented in the following
topics: Preferential Sampling; Auxiliary Information; Models for
$[X|S]$; Utility Functions; Sequential Design; and Approximations.

\section{Preferential Sampling}

Preferential sampling plays an important role in surveys routinely
carried out by Official Statistics agencies, as those developed in the
Brazilian Institute of Geography and Statistics (IBGE), the institution
that one of us is affiliated to. So we are well aware of the relevance
of this issue.
In addition to the areas of application of preferential sampling cited
by the discussants, it is important to mention the very topical area of
publication bias \citep[see][]{bayarri1993,franco2014publication}.

The link between preferential sampling and the methodology of
survey-sampling presented by discussants is also very helpful in
clarifying the similarities of approaches. We agree with Cressie and
Chambers (hereafter CC) that papers from the latter may bring important
aspects of sampling design to the context of Geostatistics. However, an
important distinction between the approaches is needed. While in the
context of survey methods the population size is generally fixed at a
finite $N$, this feature is not true in the context of Geostatistics.
In fact, in Geostatistics a finite value of $N$ may only be associated
with the discretization of a continuous process. Thus, part of the
similarity between the approaches stems from the current limitation of
many approaches to handle inference and prediction in Geostatistics
appropriately due to the use of discrete approximations.

We agree with Chipeta and Diggle (hereafter CD) that preferential
sampling is a method of adaptive design, which may depend on the
previous design without relying on the underlying process S. Similarly,
inference may be simplified after assuming that
the process $X$ is governed by the values of a spatially distributed
covariate $W$, thus rendering conditional independence between $X$ and
$S$ given $W$. The main difficulty is finding and quantifying such a
covariate. We will return to this issue in the next section.

\section{Auxiliary Information}

All discussants brought in the issue of replacement of some of the
latent and unobserved random processes by observed quantities. These
processes are a mere artifact to best represent our lack of knowledge
about the true mechanism underlying the observation processes, if such
truth exists. If we knew what variables cause our observational
processes to behave as they do and were able to measure these
quantities, we should use them.
Otherwise, models that incorporate relevant features such as smoothness
of these processes are a useful alternative.
This issue is not only relevant to Spatial Statistics but to any other
area of Statistics where absence of complete knowledge of the sources
of variation forces the use of qualitative information in the form of
latent processes, and Gaussian processes are one of the most adequate
choices as a first step.

We agree that covariates could (and should) be part of all model
components. Covariates may also be used to construct a deterministic
intensity function of the point process $X$, or as proxy of a random
process S in a model with a random intensity function. They may reduce
remaining spatial heterogeneity in the mean or covariance structure of
$S$, when available. CC suggested the use of other covariates $Z$,
assumed to be highly correlated with the process $S$, as a proxy.
This is a practical solution that needs to be used carefully. Despite
simplifications produced in inference, it is important to bear in mind
that this solution may also introduce another source of error in the
model.\looseness=-1

However, despite being a natural choice, we do not always have
available covariates. In other situations, the covariates are available
only in part of the region of interest. In these cases, it is often
necessary to interpolate values before including them in the model,
adding more uncertainty to the results.

\section{Models for $[X|S]$}

We agree that an appropriate model is crucial and robustness
considerations are even more important here.
Many discussants of \citet{diggle2010geostatistical} seemed to agree
that, in practical situations, it is unlikely that the design is
governed by a log-Gaussian Cox process. This issue is in line with our
views expressed in the first paragraph of the previous section. In some
cases, it may be enough to assume that the intensity of the point
process is somehow proportional to the underlying process at an
appropriate scale. In these cases, a log-linear intensity function may be
seen as the first option for approximation and may thus be a good
starting point to mitigate and understand the consequences of a
preferential sample.

In addition to methods for obtaining a robust design against selection
bias, alternative specifications for $[X | S]$ may also contribute to a
satisfactory model-based inference. All discussants have correctly
expressed concern about the sensitivity of the inference to the choice
of the model for the intensity function. We concur with their concern
and particularly welcome the use of non-parametric specifications to
replace the parametric specification of our paper. This is bound to
lead to a more flexible and hence robust alternative to the global
linearity imposed by the $\alpha+ \beta S$ predictor over the entire
region of interest.

One approach following this path allows the regression coefficients'
processes $\alpha$ and $\beta$ to vary locally over space. Inference
for these models may be performed in a simplified fashion, via
discretization of coefficients over sub-regions \citep
{junior2014point}, or exactly, by retaining the continuous variation of
the infinite dimensional regression coefficient process \citep
{goncalvez2015}. Another non-parametric approach was proposed by \citet
{kottas2007bayesian}, based on mixtures of Dirichlet process priors.
Their idea could be adapted to a prior for the intensity of a
log-Gaussian Cox process that (instead of being concentrated on) is
only centered at $\alpha+ \beta S$, again allowing more flexibility
for the intensity.

Zidek also suggested an interesting situation where the generating
process for the locations is based on another process $S^*$. In this
case, it is first necessary to assess the degree of dependence between
these two processes, that is, to evaluate if $[S, S^*] = [S] [S^*]$. If
so, the researcher can consider the points $X$ as ancillary for $S$ and
proceed with the non-preferential, standard inference as usual.
Otherwise, it may be necessary to establish some form of dependence for
$[S | S^*]$ to be able to proceed with inference. A~bivariate
specification for $[S, S^*] $ may be one way forward \citep
[see][]{gamerman2007dynamic, crainiceanu2008bivariate}.

In conclusion, our model for $[X|S]$ seems to be able to highlight the
effects of preferential sampling in the absence of a better
understanding of the true causes of variation or of relevant
covariates, but more flexible forms are needed.

\section{Utility Functions}

We agree with CC's suggestions about notation to improve the paper
understanding. In particular, we consider pertinent the suggestion of
making $s_d$ explicit in the utility function.

The choice of a particular utility function is a crucial step to
obtaining the optimal design. For this reason, we completely agree with
CC that the questions {\it how much?} and {\it why?} cannot be set
aside while the researcher is planning a new location sample.

Utility functions based on predictive variance reductions have long
been recognized as appropriate to measure improvements in predictive
accuracy. Zidek noted that more general evaluation (of
predictive/estimation performance) may require the use of other utility
functions. We already presented at least one alternative formulation
early in the paper to emphasize our adherence to this point and to
detach ourselves from the initial approaches based only on a single
function, evaluating predictive errors. Our methods apply equally well
to any quantifiable utility function.

The combination of different goals --- e.g., reducing predictive
errors, reducing uncertainty with respect to $S$, identifying
thresholds, reducing estimation errors and evaluating costs involved
--- allows the researcher to obtain designs in complex situations.
Examples of complex utility functions with competing goals can be found
in the recent work of \citet{muller2004optimal}, \citet
{ruiz2012evolutionary} and \citet{ferreira2015inhomogeneous} to name a few.
Note that these goals may be competing, e.g., reduction of the nugget
effect estimation error assigns more utility to regions close to
locations already sampled whereas reduction of predictive error assigns
more utility for regions far from the
points already sampled. The additional challenge in this case is to
weigh the various objectives involved.

Alternatively, the use of entropy in the utility function is usually
recommended for situations where multiple goals are involved \citep
{caselton1984optimal} but this is not an easily implementable solution
in practice.
We do believe that in practical situations the utility function needs
to be specified by whoever is in charge of the analysis, be it a
regulatory body, a government institution or a researcher, with the
help of a statistician. The responsible entity must be able to value
the worth of the information each component of the utility provides. If
one can answer, for example, why predictive variance reduction is
relevant, one must be able to assign its monetary value; otherwise, one
must reassess the relevance of each component included in the utility function.
This value is more easily combined with readily quantifiable monetary
components such as cost associated with each new location. An economist
specialized in the area of the study may be a useful addition to the
team setting up this enterprise at this stage.

\section{Sequential Design}

The approach we used for design \citep{muller1999simulation} allows the
planning of $m$ new sampling sites according to the utility function
defined by the researcher. This design plan can be sequential or not,
while the underlying process $S$ can be static or dynamic in time. In
cases where $S$ is a dynamic process, planning of a sequential design
scheme is a natural choice. However, in the case where $S$ is fixed in
time, one can also plan a sequential design, although it is not
possible to ensure that the resulting design will be optimal. Actually,
a sequential sample design can be a simple and cheap solution,
especially in situations where the costs are not fixed and can increase
before obtaining the desired new sample locations.

We anticipate difficulties to assigning a distribution for $[d|S,{\bf
x}]$ without reducing the flexibility of the model, when a
sequential-sampling-design strategy to update $\theta$ and $S$ through
$[d| {\bf x,y}]$ is built. It would be necessary to assess whether the
initial sample is preferential in cases suggested by CC with a pilot
study. If it is possible to assume that the pilot phase is
non-preferential, then this information will be ancillary to inference.
On the other hand, if the pilot sample is preferential, then it is
possible to incorporate this information by performing modifications in
$[X | S]$. Note
that this may remove the Poissonity of the model for $[X|S]$, due to
repulsions that may occur around previously selected locations. This
will imply more difficulty for likelihood approaches and opens up for
interesting research questions.

The case of a dynamic underlying process is more complex.
We visualize the following generalization of the 2-stage set-up
proposed by CD. In cases where samples are taken at different times in
a multi-stage scenario, we may have
\[
[Y,X,S]=\prod_t [Y_t|Y_{1:t-1}, X_{1:t},S_t][X_t|X_{1:t-1},
S_{1:t}][S_t|S_{1:t-1}],
\]
where $Z_{a:b} = ( Z_a , Z_{a+1} , \dots, Z_{b-1} , Z_b )$, for $a <
b$ integers, and dependence on hyperparameters was removed from the
notation as in CD. Conditional independence between observations given
the corresponding underlying processes may be assumed in some cases
leading to $ [Y_t|Y_{1:t-1}, X_{1:t},S_t] = [Y_t|X_{t},S_t]$.
Similarly, $[X_t|X_{1:t-1}, S_{1:t}] = [X_t|S_{t}]$ may be assumed for
some sampling schemes, although the general formulation may be required
in some cases (see the paragraph above).
Further simplification such as those suggested by CD may be assumed,
depending on the situation.
\citet{ferreira2015inhomogeneous} worked on a similar structure,
simplified by his non-preferential sampling scheme.

In the practical situation presented by CD, where the sampling order is
a crucial factor, a simple utility function based on predictive
variance reductions or exceedances probably would not be enough to
produce a satisfactory sample design. In challenging situations like
this, it would be necessary to choose more complex forms to reward each
sample unit, at each time, in a sequential sampling scheme.

\section{Approximations}

Questioning the stationarity assumption is a mandatory task in any
Spatial Statistics problem.
Stationarity can always be seen as an approximation to a more complex
underlying dependence. But it turns out that in many practical
situations it has proved to be a reasonable, viable solution.
Alternatively, approaches based on convolution process \citep
{higdon2002space} or the use of more flexible (non-parametric)
structures to enhance the spatial dependence structure of $S$ can be
used. Again, covariates are always relevant options to handle the
large-scale heterogeneity considered by CD.
Obviously, an increase in the complexity of the model can also
complicate the evaluation of the utility function used to obtain the
optimal design and parsimony may have to be called into action.

We agree with Zidek that other approximation of $V(S|{\bf x,y})$ could
be used in this step. A simpler, but expensive, alternative is to
generate sub-chains in order to estimate this quantity during MCMC.
Alternative analytical approaches may be preferred to avoid the
increasing computational cost. It is important to emphasize that
\begin{itemize}

\item The approximation was only used to evaluate the utility function,
since an analytical expression for $V(S|{\bf x,y})$ is not available;
the values of $S$ were always sampled from the correct distribution
$[S|{\bf y, x}]$;

\item Other utility functions, as those used in our case study, may not
require any approximation.

\end{itemize}



%
